\documentclass[10pt]{article}

\usepackage[margin=1in]{geometry}  
\usepackage{graphicx}              
\usepackage{amsmath}               
\usepackage{amsfonts}              
\usepackage{amsthm}                
\usepackage{amssymb}
\usepackage{booktabs}
\usepackage{longtable}

\DeclareRobustCommand{\divby}{%
  \mathrel{\vbox{\baselineskip.65ex\lineskiplimit0pt\hbox{.}\hbox{.}\hbox{.}}}%
}

\DeclareMathOperator{\lcm}{lcm}

\begin{document}

\nocite{*}

\title{Weighted complete intersection del Pezzo surfaces}

\author{Evgeny Mayanskiy}

\maketitle

\begin{abstract}
  We classify codimension $2$ well-formed and quasi-smooth weighted complete intersection del Pezzo surfaces.
\end{abstract}

\section{Introduction}

According to \cite{Chen3}, Theorem~$1$.$3$, a weighted complete intersection del Pezzo surface, which is well-formed, quasi-smooth and is not an intersection with a linear cone, has codimension $0$, $1$ or $2$. Weighted $2$-dimensional del Pezzo hypersurfaces were studied in \cite{JohnsonKollar}, \cite{Boyer}, \cite{CheltsovShramov} and others. Weighted codimension $2$ complete intersection del Pezzo surfaces of amplitude $1$ were classified in \cite{KimPark}.\\

The goal of this note is to classify weighted complete intersection del Pezzo surfaces of codimension $2$.\\

Our main result is the following. The reader may refer to the next section for the notation and terminology.\\

{\bf Main Theorem.} {\it Let $a_0 \leq a_1 \leq a_2 \leq a_3 \leq a_4$, $d_1\leq d_2$, $X\subset \mathbb P (a_0, a_1, a_2, a_3, a_4)$ be a well-formed quasi-smooth weighted complete intersection del Pezzo surface, given by the intersection of two quasi-homogeneous polynomials of degrees $d_1$ and $d_2$, which is not an intersection with a linear cone. Then one of the following holds:
\begin{enumerate}
\item[(1)] either $a_4 \leq 500$, $d_2\leq 1000$,
\item[(2)] or $(a_0, a_1, a_2, a_3, a_4 ; \; d_1,d_2)$ appears in Table~\ref{table:1} (see Appendix).
\end{enumerate}}

\vspace{3ex}

In other words, we explicitly classify infinite series del Pezzo surface complete intersections and leave the computation of the ``sporadic cases'' in $(1)$ to a computer. The latter is straightforward using Theorem~$WF$ and Theorem~$QS$ of Iano-Fletcher quoted in the next section. See Table~\ref{table:2} in Appendix.\\

\section{Notation and assumptions}

We work over a fixed algebraically closed field of characteristic $0$. Our terminology and notation for weighted complete intersections and projective spaces follow \cite{Fletcher}.\\

$\mathbb P=\mathbb P (a_0,\ldots , a_n)$, $a_0 \leq \ldots \leq a_n$, denotes the weighted projective space with weights $a_0,\ldots , a_n$ on the coordinates $x_0,\ldots , x_n$.\\

A weighted complete intersection $X\subset \mathbb P$ of codimension $c$ is assumed to be given by {\it general} quasi-homogeneous polynomials $F_1,\ldots , F_c$ of degrees $d_1,\ldots , d_c$, $d_1\leq \ldots \leq d_c$. The generality hypothesis simplifies the exposition and is not essential otherwise. It implies that a monomial $x_0^{m_0}\ldots x_n^{m_n}$ appears in $F_i$ with a nonzero coefficient if and only if $d_i=\sum_{j} m_j \cdot a_j$.\\

We restrict our attention only to those weighted complete intersections, which are {\it not} {\it intersections} {\it with} {\it linear} {\it cones}, i.e. no degree $d_i$ is equal to one of the weights $a_0,\ldots , a_n$.\\

Our classification is based on the following criteria proved in \cite{Fletcher}. The reader may refer to \cite{Fletcher} for the notions of well-formedness and quasi-smoothness.\\

Denote $b_{i_1\cdots i_k} = \gcd (a_0\ldots \hat{a}_{i_1} \ldots \hat{a}_{i_k} \ldots a_n)$, $i_1 < \ldots < i_k$. Let $X\subset \mathbb P (a_0,a_1,a_2,a_3,a_4)$ be a codimension $c=2$ weighted complete intersection (general and not an intersection with a linear cone).\\

{\bf Theorem~$WF$ (cf. \cite{Fletcher}, $6.11$).} {\it $X$ is well-formed if and only if 
\begin{enumerate}
\item[(1)] $\forall \; i<j<k$, $\; b_{ijk} \mid d_1$ or $\; b_{ijk} \mid d_2$,
\item[(2)] $\forall \; i<j$, $\; b_{ij} \mid d_1$ and $\; b_{ij} \mid d_2$,
\item[(3)] $\forall i$, $\; b_{i} =1$.
\end{enumerate}}

\vspace{2ex}

{\bf Theorem~$QS$ (cf. \cite{Fletcher}, $8.7$).} {\it $X$ is quasi-smooth if and only if 
\begin{enumerate}
\item[(1)] $\forall i\in \{0,1,2,3,4\}$, 
\begin{itemize}
\item either $\; d_1 \divby a_i$,
\item or $\; d_2 \divby a_i$,
\item or $\; d_1-a_e,\; d_2-a_f \in (a_i)$ for some $e\neq f$;
\end{itemize}
\item[(2)] $\forall \; 0\leq i<j \leq 4$,
\begin{itemize}
\item either $\; d_1, d_2 \in (a_i,a_j)$,
\item or $\; d_1,\; d_2-a_e \in (a_i,a_j)$ for some $e$, 
\item or $\; d_1-a_e, \; d_2 \in (a_i,a_j)$ for some $e$, 
\item or $\; d_1-a_e, \; d_2-a_f\in (a_i,a_j),\;\; \#\{e\}=\#\{f\}=2,\; \{e,f\}=\{k,l,m\}$;
\end{itemize}
\item[(3)] $\forall \; 0\leq k<l<m \leq 4$,
\begin{itemize}
\item either $\; d_1, d_2 \in (a_k,a_l, a_m)$,
\item or $\; d_1,\; d_2-a_i,\; d_2-a_j \in (a_k,a_l, a_m)$, 
\item or $\; d_1-a_i,\; d_1-a_j,\; d_2 \in (a_k,a_l, a_m)$.
\end{itemize}
\end{enumerate}}

Here we assume that $\{ i,j,k,l,m \}=\{0,1,2,3,4\}$ and denote $(a_{i_1},\ldots , a_{i_p})=\{ \sum_{s=1}^p \; u_s\cdot a_{i_s}\; \mid \; u_s\in \mathbb Z ,\; u_s\geq 0 \}$.\\

The {\it amplitude} of $X$ is defined as 
$$
I = a_0+a_1+a_2+a_3+a_4-d_1-d_2.
$$

In accordance with the adjunction formula (\cite{Dolgachev}, Theorem~$3.3.4$, \cite{Fletcher}, $6.14$), we give the following definition.\\

{\bf Definition.} A well-formed and quasi-smooth weighted complete intersection $X\subset \mathbb P (a_0,a_1,a_2,a_3,a_4)$ of dimension $2$ is a {\it del Pezzo surface}, if 
\begin{equation}\label{dPcondition}
I\geq 1.
\end{equation}

\vspace{2ex}

Throughout the note, $\lambda$, $\mu$, $\nu$, $\alpha$, $\beta$, $\gamma$ denote integers.\\

\section{Preliminary analysis}

The proof of the Main Theorem is an elementary, but somewhat lengthy analysis of all sequences $(a_0$, $a_1$, $a_2$, $a_3$, $a_4$; $d_1$, $d_2)$, which satisfy the conditions of Theorem~$WF$ and Theorem~$QS$ of Iano-Fletcher quoted above, as well as inequality (\ref{dPcondition}). In this section we list all possibilities for the expressions of $d_1$ and $d_2$ in terms of the $a_i$'s. Each of the resulting cases is studied separately in this and the next sections. All together they lead to the proof of the Main Theorem.\\

{\bf Lemma~$1$.} {\it Let $a_0 \leq a_1 \leq a_2 \leq a_3 \leq a_4$, $d_1\leq d_2$, $X\subset \mathbb P (a_0, a_1, a_2, a_3, a_4)$ be a well-formed quasi-smooth weighted complete intersection del Pezzo surface, given by the intersection of two quasi-homogeneous polynomials of degrees $d_1$ and $d_2$, which is not an intersection with a linear cone.\\

Then $(d_1,d_2)$ is one of the following:
\begin{enumerate}
\item[(01)] $d_1=a_0+a_4$, $d_2=a_1+a_4$; \qquad \qquad {\it (02)} $\; d_1=a_0+a_4$, $d_2=a_2+a_4$;
\item[(12)] $d_1=a_1+a_4$, $d_2=a_2+a_4$; \qquad \qquad {\it (03)} $\; d_1=a_0+a_4$, $d_2=a_3+a_4$;
\item[(13)] $d_1=a_1+a_4$, $d_2=a_3+a_4$; \qquad \qquad {\it (23)} $\; d_1=a_2+a_4$, $d_2=a_3+a_4$;
\item[(4.1)] $d_1=a_0+a_3$, $d_2=2a_4$;  \qquad {\it (4.2)} $\; d_1=a_1+a_3$, $d_2=2a_4$; \qquad {\it (4.3)} $\; d_1=a_2+a_3$, $d_2=2a_4$;
\item[(4.4)] $d_1=a_0+a_4$, $d_2=2a_4$;  \qquad {\it (4.5)} $\; d_1=a_1+a_4$, $d_2=2a_4$;  \qquad {\it (4.6)} $\; d_1=a_2+a_4$, $d_2=2a_4$;
\item[(4.7)] $d_1=2a_3$, $d_2=2a_4$;  \qquad \quad \; {\it (4.8)} $\; d_1=a_3+a_4$, $d_2=2a_4$; \qquad {\it (4.9)} $\; d_1=2a_4$, $d_2=2a_4$.
\end{enumerate}}

\vspace{3ex}

{\it Proof:} First, note that Theorem~$QS$~$(2)$, $(i,j)=(3,4)$, implies that $d_1\geq a_0+a_3$.\\

Hence (\ref{dPcondition}) gives 
\begin{equation}\label{dPcondition2}
d_2 < a_1+a_2+a_4.
\end{equation}

By Theorem~$QS$~$(1)$, $i=4$, we have:
\begin{itemize}
\item either $d_2 \divby a_4$, i.e. $d_2=2a_4$ by (\ref{dPcondition2}), 
\item or $d_1 \divby a_4$, i.e. $d_1=d_2=2a_4$ by (\ref{dPcondition}) and Theorem~$QS$~$(2)$, $(i,j)=(3,4)$, 
\item or $d_1=a_e+\lambda a_4,\; d_2=a_f+\mu a_4 $ for some $\; e\neq f,\; \mu, \lambda \geq 1$.
\end{itemize}

In the last case, (\ref{dPcondition}) implies that $\lambda +\mu \leq 2$, i.e. $\lambda=\mu=1$. This gives items $(01)$, $(02)$, $(12)$, $(03)$, $(13)$, $(23)$.\\

If $d_2=2a_4$, then Theorem~$QS$~$(2)$, $(i,j)=(3,4)$, requires $d_1-a_e\in (a_3,a_4)$, which gives items $(4.1)-(4.9)$. {\it QED}\\

{\bf Corollary~$1$.} {\it Under the assumptions of Lemma~$1$, $d_2\leq 2a_4$.}\\

The next lemma will be used frequently.\\

{\bf Lemma~$2$.} {\it Under the assumptions of Lemma~$1$, if $\; a_0={\lambda}_0 a_i$, $a_1={\lambda}_1 a_i$, $a_2={\lambda}_2 a_i$, $a_3={\lambda}_3 a_i$, $a_4={\lambda}_4 a_i$ for some $i$, then 
$$
a_4\leq {\lambda}_4 \cdot \min_j \lcm (\mbox{denominators of } {\lambda}_0,\ldots , \hat{\lambda}_j ,\ldots , {\lambda}_4).
$$}

{\it Proof:} By Theorem~$WF$, $a_i= \lcm (\mbox{denominators of } {\lambda}_0,\ldots , \hat{\lambda}_j ,\ldots , {\lambda}_4)$ for any $j$, if ${\lambda}_k$'s are reduced fractions. {\it QED}\\

\vspace{3ex}

{\bf Lemma~$3$.} {\it Under the assumptions of Lemma~$1$, $a_4\leq 3$, whenever three of the weights $a_i$ coincide.}\\

{\it Proof:} There are three cases.\\

{\it Case $1$.} Assume $a_0=a_1=a_2$. Then by Theorem~$WF$, $d_1 \divby a_0$ and $d_2 \divby a_0$.\\

If $d_2\geq a_3+a_4$, then $d_1=2a_0\geq a_0+a_3$, i.e. $a_0=a_1=a_2=a_3$, which should be $1$ by Theorem~$WF$. Then (\ref{dPcondition}) implies that $a_4=1$ as well.\\

Otherwise, by Lemma~$1$, $d_1=d_2=a_0+a_4$, and so $a_4 \divby a_0$. Hence by Theorem~$WF$, $a_0=a_1=a_2=1$, $d_1=d_2=1+a_4$. By (\ref{dPcondition}), $a_3=a_4$, which should be $1$ by Theorem~$WF$.\\

{\it Case $2$.} Assume $a_1=a_2=a_3$. Then by Theorem~$WF$, $d_1 \divby a_1$ and $d_2 \divby a_1$. Hence by Lemma~$1$, $2a_4 \divby a_1$, i.e. $a_1\in \{ 1,2 \}$ by Theorem~$WF$. Then (\ref{dPcondition}) and Lemma~$1$ imply that $a_4 < 2a_1 \leq 4$. \\

{\it Case $3$.} Assume $a_2=a_3=a_4$. Then $d_1=d_2=2a_4$ and (\ref{dPcondition}) gives:
$$
a_4 < a_0+a_1.
$$

By Theorem~$QS$~$(1)$, $i=1$, $2a_4 \divby a_1$. Then $a_1=2$ by Theorem~$WF$, and so $a_4 < 2a_1=4$. {\it QED}\\

\vspace{3ex}

{\bf Lemma~$4$.} {\it The} Main Theorem {\it holds, if $a_i=a_j$ for some $i\neq j$.}\\

{\it Proof:} There are four cases. In each case we consider separately the subcases arising from Lemma~$1$.\\

{\it Case $(a_0=a_1)$}: {\it Subcase $(01)$:} $d_1=d_2=a_0+a_4$, $a_4<a_2+a_3$ by (\ref{dPcondition}).\\

By Theorem~$QS$~$(1)$, $i=3$,
\begin{itemize}
\item either $a_0+a_4=\mu a_3$, where $\mu = 2$ by (\ref{dPcondition}), i.e. $a_4=2a_3-a_0$, 
\item or $a_4 \divby a_3$, and so $a_4=a_3$ by (\ref{dPcondition}).
\end{itemize}

If $a_4=a_3$, then $d_1=d_2=a_0+a_3 \divby a_3$ by Theorem~$WF$. This situation was considered in Lemma~$3$.\\

So, we may assume that $a_4=2a_3-a_0$, $d_1=d_2=2a_3$ and $a_3<a_0+a_2$. Then Theorem~$QS$~$(1)$, $i=2$, gives:
\begin{itemize}
\item either $2a_3 \divby a_2$, i.e. $a_3=a_2$ or $a_3=3a_2/2$ by (\ref{dPcondition}), 
\item or $a_3=a_2+a_0/2$.
\end{itemize}

If $a_3=a_2$, then $d_1=d_2=2a_2 \divby a_0$ by Theorem~$WF$, i.e.
\begin{itemize}
\item either $a_0=a_1=1$, $a_2=a_3=t$, $a_4=2t-1$, $d_1=d_2=2t$,
\item or $a_0=a_1=2$, $a_2=a_3=t$, $a_4=2t-2$, $d_1=d_2=2t$.
\end{itemize}

Note that in the second case $t$ must be odd by Theorem~$WF$. Both these cases appear in Table~\ref{table:1} (No.~$15$ and No.~$21$).\\

If $a_3=3a_2/2$, then $a_2<2a_0$ by (\ref{dPcondition}) and $d_1=d_2=3a_2 \divby a_0$ by Theorem~$WF$, i.e. $a_2=\lambda a_0/3$, $a_1=a_0$, $a_4=(\lambda-1)a_0$, $\lambda < 6$ by (\ref{dPcondition}), and so $a_4\leq 3(\lambda -1) < 15$ by Lemma~$2$.\\

If $a_3=a_2+a_0/2$, then $d_1=d_2=2a_2+a_0$. By Theorem~$WF$, $2a_2 \divby a_0$, and so 
$$
a_0=a_1=2,\; a_2=t, \; a_3=t+1,\; a_4=2t,\; d_1=d_2=2t+2.
$$

Since either $a_2$ or $a_3$ is even, this situation violates conditions of Theorem~$WF$.\\

{\it Subcase $(02)=(12)$:} $d_1=a_0+a_4$, $d_2=a_2+a_4$, $a_4<a_0+a_3$ by (\ref{dPcondition}).\\

By Theorem~$QS$~$(1)$, $i=3$,
\begin{itemize}
\item either $a_0+a_4 \divby a_3$, i.e. $a_4=2a_3-a_0$ by (\ref{dPcondition}),
\item or $a_2+a_4 \divby a_3$, i.e. $a_4=2a_3-a_2$ by (\ref{dPcondition}),
\item or $a_0+a_4 =a_e +\lambda a_3$ and $a_2+a_4 =a_f +\mu a_3$ for some $e\neq f$, and so $a_4=a_3$ or ($a_3=a_2$, $a_4=2a_2-a_0$) or ($a_3=2(a_2-a_0)$, $a_4=3(a_2-a_0)$).
\end{itemize}

If $a_4=a_3$, then $a_i=1$ for all $i$ by Theorem~$WF$.\\

If $a_2=a_3$, $a_4=2a_2-a_0$, then $d_1=2a_2$, $d_2=3a_2-a_0$, $a_2<2a_0$. By Theorem~$WF$, $6a_2 \divby a_0$, i.e. $a_2=a_3=\lambda a_0/6$, $\lambda < 12$, and so $a_4=(\lambda /3 -1)a_0<24$ by Lemma~$2$.\\

If $a_3=2(a_2-a_0)$, $a_4=3(a_2-a_0)$, then $d_1=a_0+3(a_2-a_0)$, $d_2=a_0+4(a_2-a_0)$. Hence $a_0 \divby a_2-a_0$ by Theorem~$WF$, i.e. 
$$
a_0=a_1=\frac{\lambda}{1+\lambda} a_2,\; \mbox{and so} \quad a_4=\frac{3}{1+\lambda} a_2\leq 3 \quad \mbox{by Lemma}\; 2.
$$

If $a_4=2a_3-a_0$, then $d_1=2a_3$, $d_2=2a_3+a_2-a_0$, $a_3<2a_0$.\\

By Theorem~$QS$~$(1)$, $i=2$,
\begin{itemize}
\item either $2a_3 \divby a_2$, i.e. $a_3=a_2$ or $a_3=3a_2/2$ by (\ref{dPcondition}),
\item or $a_3=a_2+a_0/2$.
\end{itemize}

The case $a_2=a_3$ was considered above. If $a_3=3a_2/2$, then $d_1=3a_2$, $d_2=4a_2-a_0$, $a_2<4a_0/3$. By Theorem~$WF$, $12a_2\divby a_0$, i.e. $a_2=\lambda a_0/12$, $\lambda < 16$, and so $a_4=(\lambda /4-1) a_0 <48$ by Lemma~$2$.\\

If $a_3=a_2+a_0/2$, then $d_1=2a_2+a_0$, $d_2=3a_2$, $a_2<3a_0/2$. By Theorem~$WF$, $6a_2\divby a_0$, i.e. $a_2=\lambda a_0/6$, $\lambda < 9$, and so $a_4=(\lambda /3) a_0 <18$ by Lemma~$2$.\\

If $a_4=2a_3-a_2$, then $d_1=2a_3-a_2+a_0$, $d_2=2a_3$, $a_3<a_0+a_2$.\\

By Theorem~$QS$~$(1)$, $i=2$,
\begin{itemize}
\item either $2a_3 \divby a_2$, i.e. $a_3=a_2$ or $a_3=3a_2/2$ by (\ref{dPcondition}),
\item or $2a_3+a_0 \divby a_2$, i.e. $a_3=2a_2-a_0/2$ or $a_3=3a_2/2-a_0/2$ by (\ref{dPcondition}),
\item or $a_3=2a_2-a_0$.
\end{itemize}

The case $a_2=a_3=a_4$ was considered in Lemma~$3$.\\

If $a_3=3a_2/2$, then $d_1=2a_2+a_0$, $d_2=3a_2$, $a_2<2a_0$. By Theorem~$WF$, $6a_2\divby a_0$, i.e. $a_2=\lambda a_0/6$, $\lambda < 12$, and so $a_4=(\lambda /3) a_0 <24$ by Lemma~$2$.\\

If $a_3=2a_2-a_0/2$, then $d_1=3a_2$, $d_2=4a_2-a_0$, $a_2<3a_0/2$. By Theorem~$WF$, $12a_2\divby a_0$, i.e. $a_2=\lambda a_0/12$, $\lambda < 18$, and so $a_4=(\lambda /4-1) a_0 <54$ by Lemma~$2$.\\

If $a_3=3a_2/2-a_0/2$, then $d_1=2a_2$, $d_2=3a_2-a_0$, $a_2<3a_0$. By Theorem~$WF$, $6a_2\divby a_0$, i.e. $a_2=\lambda a_0/6$, $\lambda < 18$, and so $a_4=(\lambda /3-1) a_0 <36$ by Lemma~$2$.\\

If $a_3=2a_2-a_0$, then $d_1=3a_2-a_0$, $d_2=4a_2-2a_0$, $a_2<2a_0$. By Theorem~$WF$, $12a_2\divby a_0$, i.e. $a_2=\lambda a_0/12$, $\lambda < 24$, and so $a_4=(\lambda /4-2) a_0 <72$ by Lemma~$2$.\\

{\it Subcases $(03)=(13)$, $(23)$, $(4.1)-(4.9)$:} The same analysis as in the previous subcases gives:
$$
a_4<60.
$$

{\it Case $(a_1=a_2)$}: {\it Subcase $(01)=(02)$:} $d_1=a_0+a_4$, $d_2=a_1+a_4$, $a_4<a_1+a_3$ by (\ref{dPcondition}).\\

By Theorem~$QS$~$(1)$, $i=3$,
\begin{itemize}
\item either $a_0+a_4 \divby a_3$, i.e. $a_4=2a_3-a_0$ by (\ref{dPcondition}),
\item or $a_1+a_4 \divby a_3$, i.e. $a_4=2a_3-a_1$ by (\ref{dPcondition}),
\item or $a_0+a_4 =a_e +\lambda a_3$ and $a_1+a_4 =a_f +\mu a_3$ for some $e\neq f$, and so $a_4=a_3$ or ($a_3=2(a_1-a_0)$, $a_4=3(a_1-a_0)$).
\end{itemize}

If $a_3=2(a_1-a_0)$, $a_4=3(a_1-a_0)$, then $d_1=a_0+3(a_1-a_0)$, $d_2=a_0+4(a_1-a_0)$. Hence $a_0\divby a_1-a_0$ by Theorem~$WF$, and so 
$$
a_1=a_2=\frac{1+\lambda}{\lambda}a_0,\; a_4=\frac{3}{\lambda}a_0\leq 3 \quad \mbox{by Lemma}\; 2.
$$

If $a_3=a_4$, then $a_1=a_4$ by Theorem~$WF$. This situation was considered in Lemma~$3$.\\

If $a_4=2a_3-a_1$, then $d_1=2a_3+a_0-a_1$, $d_2=2a_3$, $a_3<2a_1$. By Theorem~$WF$,
\begin{itemize}
\item either $2a_3 \divby a_1$, i.e. $a_3=3a_1/2$ by (\ref{dPcondition}),
\item or $2a_3+a_0 \divby a_1$, i.e. $a_3=2a_1-a_0/2$ or $a_3=3a_1/2-a_0/2$ by (\ref{dPcondition}).
\end{itemize}

If $a_3=3a_1/2$, then $a_4=2a_1\leq 4$ by Lemma~$2$.\\

If $a_3=2a_1-a_0/2$, then $d_1=3a_1$, $d_2=4a_1-a_0$. By Theorem~$QS$~$(1)$, $i=0$,
\begin{itemize}
\item either $3a_1 \divby a_0$,
\item or $4a_1 \divby a_0$, i.e. $a_1=a_2=\lambda a_0/4$, $a_3=(\lambda -1)a_0/2$, $a_4=(3\lambda /4-1)a_0$, $d_1=3\lambda a_0/4$, $d_2=(\lambda-1)a_0$.
\end{itemize}

If $3a_1 \divby a_0$, then $a_0$, $a_1=a_2$ and $a_4$ are even, which violates conditions of Theorem~$WF$.\\

Hence 
\begin{itemize}
\item either $a_0=2$, $a_1=a_2=t$, $a_3=2t-1$, $a_4=3t-2$, $d_1=3t$, $d_2=4t-2$, where $t$ is odd,
\item or $a_0=4$, $a_1=a_2=2t+1$, $a_3=4t$, $a_4=6t-1$, $d_1=6t+3$, $d_2=8t$.
\end{itemize}

These solutions appear in Table~\ref{table:1} (No.~$25$ and No.~$35$).\\

If $a_3=3a_1/2-a_0/2$, then $d_1=2a_1$, $d_2=3a_1-a_0$. By Theorem~$QS$~$(1)$, $i=0$,
\begin{itemize}
\item either $2a_1 \divby a_0$, i.e. $a_1=a_2=\lambda a_0/2$, $a_3=(3\lambda-2)a_0/4$, $a_4=(\lambda -1)a_0$, $d_1=\lambda a_0$, $d_2=(3\lambda /2-1)a_0$,
\item or $3a_1 \divby a_0$, i.e. $a_1=a_2=\lambda a_0/3$, $a_3=(\lambda -1)a_0/2$, $a_4=(2\lambda /3-1)a_0$, $d_1=2\lambda a_0/3$, $d_2=(\lambda-1)a_0$.
\end{itemize}

Hence 
\begin{itemize}
\item either $a_0=1$, $a_1=a_2=2t+1$, $a_3=3t+1$, $a_4=4t+1$, $d_1=4t+2$, $d_2=6t+2$,
\item or $a_0=3$, $a_1=a_2=2t+1$, $a_3=3t$, $a_4=4t-1$, $d_1=4t+2$, $d_2=6t$.
\end{itemize}

Note that in the second case $t\not\equiv 1 \pmod 3$ by Theorem~$WF$. These solutions appear in Table~\ref{table:1} (No.~$19$ and No.~$30$).\\

If $a_4=2a_3-a_0$, then $d_1=2a_3$, $d_2=2a_3+a_1-a_0$, $a_3<a_0+a_1$. By Theorem~$WF$,
\begin{itemize}
\item either $2a_3 \divby a_1$, i.e. $a_3=3a_1/2$ by (\ref{dPcondition}), 
\item or $2a_3-a_0 \divby a_1$, i.e. $a_3=a_1+a_0/2$ by (\ref{dPcondition}).
\end{itemize}

If $a_3=a_1+a_0/2$, then $a_4=2a_1$, $d_1=2a_1+a_0$, $d_2=3a_1$. By Theorem~$WF$, $a_0=a_1=a_2$. This situation was considered in Lemma~$3$.\\

If $a_3=3a_1/2$, then $d_1=3a_1$, $d_2=4a_1-a_0$, $a_1<2a_0$. By Theorem~$QS$~$(1)$, $i=0$, either $12a_1 \divby a_0$ or $5a_1 \divby a_0$. If $12a_1 \divby a_0$, then $a_1=a_2=\lambda a_0/12$, $\lambda < 24$, and so $a_4=(\lambda /4-1)a_0<60$ by Lemma~$2$. If $5a_1 \divby a_0$, then $a_1=a_2=\lambda a_0/5$, $\lambda < 10$, and so $a_4=(3\lambda /5-1)a_0<30$ by Lemma~$2$.\\

{\it Subcase $(12)$:} $d_1=d_2=a_1+a_4$, $a_4<a_0+a_3$ by (\ref{dPcondition}).\\

By Theorem~$QS$~$(1)$, $i=3$,
\begin{itemize}
\item either $a_1+a_4 \divby a_3$, i.e. $a_4=2a_3-a_1$ by (\ref{dPcondition}),
\item or $a_4 \divby a_3$, i.e. $a_4=a_3$ by (\ref{dPcondition}).
\end{itemize}

If $a_4=2a_3-a_1$, then $d_1=d_2=2a_3$, $a_3<a_0+a_1$. By Theorem~$WF$, $2a_3\divby a_1$, i.e. $a_3=3a_1/2$ by (\ref{dPcondition}). Then $a_4=2a_1\leq 4$ by Lemma~$2$.\\

If $a_3=a_4$, then $d_1=d_2=a_1+a_3$. By Theorem~$WF$, $a_3\divby a_1$, i.e. $a_3=a_4=\lambda a_1$. Hence $a_0=a_1=a_2=1$ by Lemma~$2$. This situation was considered in Lemma~$3$.\\

{\it Subcase $(03)$:} $d_1=a_0+a_4$, $d_2=a_3+a_4$, $a_4<2a_1$ by (\ref{dPcondition}).\\

By Theorem~$QS$~$(1)$, $i=3$,
\begin{itemize}
\item either $a_4 \divby a_3$, i.e. $a_4=a_3$ by (\ref{dPcondition}),
\item or $a_0 + a_4 \divby a_3$, i.e. $a_4=2a_3-a_0$ by (\ref{dPcondition}),
\item or $a_4 =a_e +\lambda a_3$ and $a_0+a_4 =a_f +\mu a_3$ for some $e\neq f$, and so $a_4=a_3+a_1-a_0$.
\end{itemize}

If $a_4=a_3+a_1-a_0$, then $d_1=a_1+a_3$, $d_2=2a_3+a_1-a_0$, $a_3<a_0+a_1$. By Theorem~$WF$,
\begin{itemize}
\item either $a_3 \divby a_1$, i.e. $a_1=a_2=a_3$ by (\ref{dPcondition}),
\item or $2a_3-a_0 \divby a_1$, i.e. $a_3=a_1+a_0/2$ by (\ref{dPcondition}).
\end{itemize}

If $a_3=a_1+a_0/2$, then $d_1=2a_1+a_0/2$, $d_2=3a_1$. By Theorem~$QS$~$(1)$, $i=0$,
\begin{itemize}
\item either $4a_1 \divby a_0$, i.e. $a_1=a_2=\lambda a_0/4$, $a_3=(\lambda +2)a_0/4$, $a_4=(\lambda -1)a_0/2$, $d_1=(\lambda+1)a_0/2$, $d_2=3\lambda a_0/4$, 
\item or $3a_1\divby a_0$.
\end{itemize}

If $3a_1\divby a_0$, then $a_0$, $a_1=a_2$ are even. Hence $d_1$ should be even by Theorem~$WF$. This implies that $4 \mid a_0$, and so $a_4$ is even too. This violates conditions of Theorem~$WF$.\\

If $4a_1 \divby a_0$ and $a_0=2$, then a condition of Theorem~$WF$ is violated. Hence 
$$
a_0=4,\; a_1=a_2=2t+1,\; a_3=2t+3,\; a_4=4t,\; d_1=4t+4,\; d_2=6t+3.
$$

This solution appears in Table~\ref{table:1} (No.~$34$).\\

If $a_4=2a_3-a_0$, then $d_1=2a_3$, $d_2=3a_3-a_0$, $a_3<a_1+a_0/2$. By Theorem~$WF$, 
\begin{itemize}
\item either $2a_3 \divby a_1$, i.e. $a_1=a_2=a_3$,
\item or $3a_3-a_0\divby a_1$, i.e. $a_3=a_1+a_0/3$.
\end{itemize}

If $a_3=a_1+a_0/3$, then $d_1=2a_1+2a_0/3$, $d_2=3a_1$. By Theorem~$QS$~$(1)$, $i=0$, $6a_1\divby a_0$, i.e. 
$$
a_1=a_2=\frac{\lambda a_0}{6},\; a_3=\frac{\lambda +2}{6}a_0,\; a_4=\frac{\lambda -1}{3}a_0,\; d_1=\frac{\lambda +2}{3}a_0,\; d_2=\frac{\lambda a_0}{2}.
$$

This violates the condition of Theorem~$QS$~$(2)$, $(i,j)=(1,2)$.\\

If $a_4=a_3$, then $d_1=a_0+a_3$, $d_2=2a_3$, $a_3<2a_1$. By Theorem~$WF$, 
\begin{itemize}
\item either $2a_3 \divby a_1$, i.e. $a_3=3a_1/2$,
\item or $a_3+a_0\divby a_1$, i.e. $a_3=2a_1-a_0$.
\end{itemize}

If $a_3=3a_1/2$, then $a_1=2$ by Theorem~$WF$, and so $a_4=3$. If $a_3=2a_1-a_0$, then $d_1=2a_1$, $d_2=4a_1-2a_0$. By Theorem~$QS$~$(1)$, $i=0$, 
\begin{itemize}
\item either $4a_1\divby a_0$, i.e. $a_1=a_2=\lambda a_0/4$, $a_3=a_4=(\lambda-2)a_0/2$, $d_1=\lambda a_0/2$, $d_2=(\lambda-2)a_0$,
\item or $3a_1\divby a_0$, i.e. $a_1=a_2=\lambda a_0/3$, $a_3=a_4=(2\lambda /3-1)a_0$, $d_1=2\lambda a_0/3$, $d_2=(4\lambda /3-2)a_0$.
\end{itemize}

If $a_0=3$ in the second case, then either the condition of Theorem~$QS$~$(2)$, $(i,j)=(1,2)$, or Theorem~$WF$ is violated.\\

If $a_0=1$ in the first case, then by Theorem~$QS$~$(2)$, $(i,j)=(1,2)$, $\lambda \in \{ 4,8,12 \}$, i.e. $a_4\leq 5$.\\

If $a_0=2$ in the first case, then by Theorem~$QS$~$(2)$, $(i,j)=(1,2)$, $\lambda =6$, i.e. $a_4\leq 4$.\\

If $a_0=4$ in the first case, then the condition of Theorem~$QS$~$(2)$, $(i,j)=(1,2)$ is violated.\\

{\it Subcases $(13)=(23)$, $(4.1)-(4.9)$:} The same analysis as in the previous subcases gives:
$$
a_4<40.
$$

Note that in Subcase $(4.1)$ it may happen that $a_3=a_4=2a_1-a_0$, $d_1=2a_1$, $d_2=4a_1-2a_0$. Then by Theorem~$QS$~$(2)$, $(i,j)=(1,2)$, either $a_1=3a_0$ or $a_1=2a_0$ or $a_1=3a_0/2$ or $a_1=a_0$. In this case, $a_4\leq 5$ by Theorem~$WF$ and Lemma~$2$.\\

{\it Case $(a_2=a_3)$}: {\it Subcase $(01)$:} $d_1=a_0+a_4$, $d_2=a_1+a_4$, $a_4<2a_2$ by (\ref{dPcondition}). By Theorem~$WF$,
\begin{itemize}
\item either $a_0+a_4\divby a_2$, i.e. $a_4=2a_2-a_0$ by (\ref{dPcondition}),
\item or $a_1+a_4\divby a_2$, i.e. $a_4=2a_2-a_1$ by (\ref{dPcondition}).
\end{itemize}

If $a_4=2a_2-a_1$, then $d_1=2a_2+a_0-a_1$, $d_2=2a_2$. By Theorem~$QS$~$(2)$, $(i,j)=(2,3)$, $a_2=a_3=2a_1-a_0$. Then $d_1=3a_1-a_0$, $d_2=4a_1-2a_0$ and by Theorem~$QS$~$(1)$, $i=0$,  
\begin{itemize}
\item either $3a_1\divby a_0$, i.e. $a_1=\lambda a_0/3$, $a_2=a_3=(2\lambda /3-1)a_0$, $a_4=(\lambda-2)a_0$, $d_1=(\lambda-1)a_0$, $d_2=(4\lambda /3 -2)a_0$,
\item or $4a_1\divby a_0$, i.e. $a_1=\lambda a_0/4$, $a_2=a_3=(\lambda /2-1)a_0$, $a_4=(3\lambda /4-2)a_0$, $d_1=(3\lambda /4-1)a_0$, $d_2=(\lambda -2)a_0$.
\end{itemize}

In the second item $a_0=1$ by Theorem~$WF$.\\

Hence
\begin{itemize}
\item either $a_0=3$, $a_1=t$, $a_2=a_3=2t-3$, $a_4=3t-6$, $d_1=3t-3$, $d_2=4t-6$, $t\not \equiv 0 \pmod 3$,
\item or $a_0=1$, $a_1=t$, $a_2=a_3=2t-1$, $a_4=3t-2$, $d_1=3t-1$, $d_2=4t-2$.
\end{itemize}

These solutions appear in Table~\ref{table:1} (No.~$17$ and No.~$26$).\\

If $a_4=2a_2-a_0$, then $d_1=2a_2$, $d_2=2a_2+a_1-a_0$. By Theorem~$QS$~$(2)$, $(i,j)=(2,3)$, $a_1=2a_0$, and so $d_1=2a_2$, $d_2=2a_2+a_0$. By Theorem~$QS$~$(1)$, $i=0$, $2a_2\divby a_0$, i.e.
$$
a_1=2a_0,\; a_2=a_3=\lambda a_0/2 ,\; a_4=(\lambda -1) a_0 ,\; d_1=\lambda a_0,\; d_2=(\lambda +1) a_0.
$$

Hence
\begin{itemize}
\item either $a_0=2$, $a_1=4$, $a_2=a_3=2t+1$, $a_4=4t$, $d_1=4t+2$, $d_2=4t+4$,
\item or $a_0=1$, $a_1=2$, $a_2=a_3=t$, $a_4=2t-1$, $d_1=2t$, $d_2=2t+1$, where $t$ is odd.
\end{itemize}

These solutions appear in Table~\ref{table:1} (No.~$16$ and No.~$24$).\\

{\it Subcase $(02)=(03)$:} $d_1=a_0+a_4$, $d_2=a_3+a_4$, $a_4<a_1+a_3$ by (\ref{dPcondition}). By Theorem~$WF$,
\begin{itemize}
\item either $a_4\divby a_3$, i.e. $a_4=a_3=a_2$ by (\ref{dPcondition}),
\item or $a_0+a_4\divby a_3$, i.e. $a_4=2a_3-a_0$ by (\ref{dPcondition}).
\end{itemize}

By Lemma~$3$, we may assume that $a_4=2a_3-a_0$, and so $d_1=2a_3$, $d_2=3a_3-a_0$, $a_3<a_0+a_1$. By Theorem~$QS$~$(1)$, $i=1$, 
\begin{itemize}
\item either $2a_3\divby a_1$, i.e. $a_3=3a_1/2$ by (\ref{dPcondition}),
\item or $3a_3-a_0\divby a_1$, i.e. $a_3=a_1+a_0/3$ or $a_3=4a_1/3+a_0/3$ by (\ref{dPcondition}),
\item or $2a_3-a_0\divby a_1$, i.e. $a_3=a_1+a_0/2$ by (\ref{dPcondition}).
\end{itemize}

If $a_3=3a_1/2$ or $a_3=4a_1/3+a_0/3$, then $a_4<90$ by Theorem~$QS$~$(1)$, $i=0$, and Lemma~$2$.\\

If $a_3=a_1+a_0/3$, then $d_1=2a_1+2a_0/3$, $d_2=3a_1$. By Theorem~$QS$~$(1)$, $i=0$, $6a_1\divby a_0$, i.e. 
$$
a_1=\frac{\lambda a_0}{6},\; a_2=a_3=\frac{\lambda +2}{6}a_0,\; a_4=\frac{\lambda -1}{3}a_0,\; d_1=\frac{\lambda +2}{3}a_0,\; d_2=\frac{\lambda a_0}{2}.
$$

Hence
\begin{itemize}
\item either $a_0=6$, $a_1=2t+1$, $a_2=a_3=2t+3$, $a_4=4t$, $d_1=4t+6$, $d_2=6t+3$, where $t\equiv 1\pmod 3$,
\item or $a_0=3$, $a_1=t$, $a_2=a_3=t+1$, $a_4=2t-1$, $d_1=2t+2$, $d_2=3t$, where $t\not \equiv -1\pmod 3$.
\end{itemize}

These solutions appear in Table~\ref{table:1} (No.~$28$ and No.~$40$).\\

If $a_3=a_1+a_0/2$, then $a_4=2a_1$, $d_1=2a_1+a_0$, $d_2=3a_1+a_0/2$. By Theorem~$WF$, $a_0=a_1$. This situation was considered in a previous case.\\

{\it Subcase $(12)=(13)$:} $d_1=a_1+a_4$, $d_2=a_3+a_4$, $a_4<a_0+a_3$ by (\ref{dPcondition}). By Theorem~$WF$,
\begin{itemize}
\item either $a_4\divby a_3$, i.e. $a_4=a_3=a_2$ by (\ref{dPcondition}),
\item or $a_1+a_4\divby a_3$, i.e. $a_4=2a_3-a_1$ by (\ref{dPcondition}).
\end{itemize}

By Lemma~$3$, we may assume that $a_4=2a_3-a_1$, and so $d_1=2a_3$, $d_2=3a_3-a_1$, $a_3<a_0+a_1$. By Theorem~$QS$~$(1)$, $i=1$, 
\begin{itemize}
\item either $2a_3\divby a_1$, i.e. $a_3=3a_1/2$ by (\ref{dPcondition}),
\item or $3a_3\divby a_1$, i.e. $a_3=4a_1/3$ or $a_3=5a_1/3$ by (\ref{dPcondition}),
\item or $3a_3-a_0\divby a_1$, i.e. $a_3=a_1+a_0/3$ or $a_3=4a_1/3+a_0/3$ by (\ref{dPcondition}).
\end{itemize}

In each case, except for $a_3=a_1+a_0/3$, $a_4<96$ by Theorem~$QS$~$(1)$, $i=0$, and Lemma~$2$.\\

If $a_3=a_1+a_0/3$, then $d_1=2a_1+2a_0/3$, $d_2=2a_1+a_0$. By Theorem~$QS$~$(1)$, $i=0$, $6a_1\divby a_0$, i.e. 
$$
a_1=\frac{\lambda a_0}{6},\; a_2=a_3=\frac{\lambda +2}{6}a_0,\; a_4=\frac{\lambda +4}{6}a_0,\; d_1=\frac{\lambda +2}{3}a_0,\; d_2=\frac{\lambda +3}{3}a_0.
$$

Hence
\begin{itemize}
\item either $a_0=6$, $a_1=2t-1$, $a_2=a_3=2t+1$, $a_4=2t+3$, $d_1=4t+2$, $d_2=4t+4$, where $t\equiv -1\pmod 3$,
\item or $a_0=3$, $a_1=t$, $a_2=a_3=t+1$, $a_4=t+2$, $d_1=2t+2$, $d_2=2t+3$, where $3\mid t$.
\end{itemize}

These solutions appear in Table~\ref{table:1} (No.~$27$ and No.~$39$).\\

{\it Subcases $(23)$, $(4.1)-(4.9)$:} The same analysis gives:
$$
a_4<30.
$$

{\it Case $(a_3=a_4)$}: {\it Subcases $(01)$, $(02)$, $(12)$:} By Theorem~$WF$, $d_1\divby a_4$ or $d_2\divby a_4$. Hence $a_2=a_3=a_4$, and we apply Lemma~$3$.\\

{\it Subcase $(03)=(4.1)=(4.4)$:} $d_1=a_0+a_4$, $d_2=2a_4$, $a_4<a_1+a_2$.  By Theorem~$QS$~$(1)$, $i=2$, 
\begin{itemize}
\item either $2a_4\divby a_2$, i.e. $a_4=3a_2/2$ by (\ref{dPcondition}),
\item or $a_0+a_4\divby a_2$, i.e. $a_4=2a_2-a_0$ by (\ref{dPcondition}),
\item or $a_0+a_4=a_e+\mu a_2$, $2a_4=a_f+\nu a_2$ for some $e\neq f$, i.e. ($a_2=2a_1-3a_0$, $a_4=3a_1-4a_0$) or ($a_1=3a_0/2$, $a_4=a_2+a_0/2$).
\end{itemize}

If $a_4=a_3=3a_2/2$, then $d_1=a_0+a_4\divby a_2/2$ by Theorem~$WF$, i.e. $a_2=2a_0$, $a_4=a_3=3a_0\leq 3$ by Lemma~$2$.\\

If $a_4=2a_2-a_0$, then $d_1=2a_2$, $d_2=4a_2-2a_0$, $a_2<a_0+a_1$. By Theorem~$QS$~$(1)$, $i=1$, 
\begin{itemize}
\item either $2a_2\divby a_1$, i.e. $a_2=3a_1/2$ by (\ref{dPcondition}),
\item or $4a_2-2a_0\divby a_1$, i.e. $a_2=\lambda a_1/4+a_0/2$, $\lambda \in \{ 3,4,5 \}$ by (\ref{dPcondition}).
\end{itemize}

If $a_2=3a_1/2$, then $d_1=3a_1$, $d_2=6a_1-2a_0$, $a_1<2a_0$. By Theorem~$QS$~$(1)$, $i=0$, either $6a_1\divby a_0$ or $5a_1\divby a_0$ or $9a_1\divby a_0$. Then $a_4<45$ by Lemma~$2$.\\ 

If $a_2=5a_1/4+a_0/2$, then $d_1=5a_1/2+a_0$, $d_2=5a_1$, $a_1<2a_0$. By Theorem~$QS$~$(1)$, $i=0$, either $15a_1\divby a_0$ or $4a_1\divby a_0$. Then $a_4 < 60$ by Lemma~$2$.\\ 

If $a_2=a_1+a_0/2$, then $a_3=a_4=2a_1$. By Theorem~$WF$, $d_1=2a_1+a_0\divby a_1$, i.e. $a_0=a_1$, which is a previously considered case.\\ 

If $a_2=3a_1/4+a_0/2$, then $a_3=a_4=3a_1/2$. By Theorem~$WF$, $d_1=3a_1/2+a_0\divby a_1/2$, i.e. $a_1=2a_0$, and so $a_4=a_3=3a_0\leq 3$ by Lemma~$2$.\\ 

If $a_2=2a_1-3a_0$, $a_4=3a_1-4a_0$, then $d_1=3(a_1-a_0)$, $d_2=6a_1-8a_0$.  By Theorem~$QS$~$(1)$, $i=1$, either $9a_0\divby a_1$ or $8a_0\divby a_1$. Then $a_4 < 27$ by Lemma~$2$.\\ 

If $a_1=3a_0/2$, $a_4=a_2+a_0/2$, then $d_1=a_2+3a_0/2$, $d_2=2a_2+a_0$.  By Theorem~$QS$~$(1)$, $i=1$, $6a_2\divby a_1$, i.e. 
$$
a_1=\frac{3a_0}{2},\; a_2=\frac{\lambda a_0}{4},\; a_3=a_4=\frac{\lambda +2}{4}a_0,\; d_1=\frac{\lambda +6}{4}a_0,\; d_2=\frac{\lambda +2}{2}a_0.
$$ 

Hence
\begin{itemize}
\item either $a_0=2$, $a_1=3$, $a_2=3t$, $a_3=a_4=3t+1$, $d_1=3t+3$, $d_2=6t+2$,
\item or $a_0=4$, $a_1=6$, $a_2=6t-3$, $a_3=a_4=6t-1$, $d_1=6t+3$, $d_2=12t-2$.
\end{itemize}

These solutions appear in Table~\ref{table:1} (No.~$22$ and No.~$32$).\\

{\it Subcase $(13)=(4.2)=(4.5)$:} $d_1=a_1+a_4$, $d_2=2a_4$, $a_4<a_0+a_2$.  By Theorem~$QS$~$(1)$, $i=2$, 
\begin{itemize}
\item either $2a_4\divby a_2$, i.e. $a_4=3a_2/2$ by (\ref{dPcondition}),
\item or $a_1+a_4\divby a_2$, i.e. $a_4=2a_2-a_1$ by (\ref{dPcondition}).
\end{itemize}

If $a_4=3a_2/2$, then  $d_1=a_1+3a_2/2$, $d_2=3a_2$, $a_2<2a_0$. By Theorem~$QS$~$(1)$, $i=1$, 
\begin{itemize}
\item either $3a_2\divby a_1$, i.e. $a_2=4a_1/3$ or $a_2=5a_1/3$ by (\ref{dPcondition}),
\item or $2a_2\divby a_1$, i.e. $a_2=3a_1/2$ by (\ref{dPcondition}),
\item or $3a_2-a_0\divby a_1$, i.e. $a_2=a_1+a_0/3$ or $a_2=4a_1/3+a_0/3$ by (\ref{dPcondition}).
\end{itemize}

In each case, $a_4<30$ by Theorem~$QS$~$(1)$, $i=0$, and Lemma~$2$.\\

If $a_4=2a_2-a_1$, then $d_1=2a_2$, $d_2=4a_2-2a_1$, $a_2<a_0+a_1$. By Theorem~$QS$~$(1)$, $i=1$, 
\begin{itemize}
\item either $4a_2\divby a_1$, i.e. $a_2=\lambda a_1/4$, $\lambda \in \{ 5,6,7 \}$ by (\ref{dPcondition}),
\item or $3a_2\divby a_1$, i.e. $a_2=\lambda a_1/3$, $\lambda \in \{ 4,5 \}$ by (\ref{dPcondition}),
\item or $4a_2-a_0\divby a_1$, i.e. $a_2=\lambda a_1/4+a_0/4$, $\lambda \in \{ 4,5,6 \}$ by (\ref{dPcondition}).
\end{itemize}

In each case, except for the last item with $\lambda =4$, $a_4<50$ by Theorem~$QS$~$(1)$, $i=0$, and Lemma~$2$.\\

If $a_2=a_1+a_0/4$, then $d_1=2a_1+a_0/2$, $d_2=2a_1+a_0$. By Theorem~$QS$~$(1)$, $i=0$, $4a_1\divby a_0$, i.e. $a_1=\lambda a_0/4$, $a_2=(\lambda +1)a_0/4$, $a_3=a_4=(\lambda +2)a_0/4$. By Theorem~$WF$,
$$
a_0=4,\; a_1=\lambda,\; a_2=\lambda +1,\; a_3=a_4=\lambda +2,\; d_1=2\lambda +2,\; d_2=2\lambda +4.
$$ 

By Theorem~$QS$~$(2)$, $(i,j)=(0,2)$, $\lambda =4t+1$. This solution appears in Table~\ref{table:1} (No.~$37$).\\

{\it Subcase $(23)=(4.3)=(4.6)$:} $d_1=a_2+a_4$, $d_2=2a_4$, $a_4<a_0+a_1$.  By Theorem~$QS$~$(1)$, $i=2$, 
\begin{itemize}
\item either $2a_4\divby a_2$, i.e. $a_4=3a_2/2$ by (\ref{dPcondition}),
\item or $2a_4-a_0\divby a_2$, i.e. $a_4=a_2+a_0/2$ by (\ref{dPcondition}),
\item or $2a_4-a_1\divby a_2$, i.e. $a_4=a_2+a_1/2$ by (\ref{dPcondition}).
\end{itemize}

If $a_4=a_2+a_0/2$, then $d_1=2a_2+a_0/2$, $d_2=2a_2+a_0$, $a_2<a_1+a_0/2$. By Theorem~$QS$~$(1)$, $i=1$, 
\begin{itemize}
\item either $2a_2+a_0/2\divby a_1$, i.e. $a_2=3a_1/2-a_0/4$ by (\ref{dPcondition}),
\item or $2a_2+a_0\divby a_1$, i.e. $a_2=3a_1/2-a_0/2$ by (\ref{dPcondition}),
\item or $2a_2+a_0/2=a_e+\mu a_1$, $2a_2+a_0=a_f+\nu a_1$ for some $e\neq f$, i.e. $a_1=5a_0/4$, $a_2=3a_0/2$ by (\ref{dPcondition}).
\end{itemize}

In each case, $a_4<16$ by Theorem~$QS$~$(1)$, $i=0$, and Lemma~$2$.\\

If $a_4=a_2+a_1/2$, then $d_1=2a_2+a_1/2$, $d_2=2a_2+a_1$, $a_2<a_0+a_1/2$. By Theorem~$QS$~$(1)$, $i=1$, $4a_2\divby a_1$, i.e. $a_2=\lambda a_1/4$, $\lambda <6$, and so $a_4=a_3=(\lambda +2)a_1/4 < 8$ by Lemma~$2$.\\

If $a_4=3a_2/2$, then $d_1=5a_2/2$, $d_2=3a_2$, $a_2<2(a_0+a_1)/3$. By Theorem~$QS$~$(1)$, $i=1$, $5a_2\divby a_1$, i.e. $a_2=6a_1/5$. Then $a_4=a_3=9a_1/5 \leq 9$ by Lemma~$2$.\\

{\it Subcase $(4.7)=(4.8)=(4.9)$:} $d_1=d_2=2a_4$, $a_4<(a_0+a_1+a_2)/2$.  By Theorem~$QS$~$(1)$, $i=2$, $2a_4\divby a_2$, i.e. $a_2=a_3=a_4$ and we apply Lemma~$3$. {\it QED}\\

\vspace{3ex}

{\bf Lemma~$5$.} {\it The} Main Theorem {\it holds in cases $(4.1)-(4.9)$ of Lemma~$1$.} \\

{\it Proof:} By Lemma~$4$, we may assume that
\begin{equation}\label{AllNonequal}
1\leq a_0 \lneqq a_1 \lneqq a_2 \lneqq a_3 \lneqq a_4.
\end{equation}

{\it Case $(4.1)$:} $d_1=a_0+a_3$, $d_2=2a_4$, $a_4<a_1+a_2$ by (\ref{dPcondition}). By Theorem~$QS$~$(1)$, $i=3$, 
\begin{itemize}
\item either $2a_4\divby a_3$, i.e. $a_4=3a_3/2$ by (\ref{dPcondition}),
\item or $2a_4-a_2\divby a_3$, i.e. $a_4=a_3+a_2/2$ by (\ref{dPcondition}),
\item or $2a_4-a_1\divby a_3$, i.e. $a_4=a_3+a_1/2$ by (\ref{dPcondition}).
\end{itemize}

If $a_4=3a_3/2$, then $d_1=a_0+a_3$, $d_2=3a_3$, $a_3<2(a_1+a_2)/3$. By Theorem~$QS$~$(1)$, $i=2$, 
\begin{itemize}
\item either $3a_3\divby a_2$, i.e. $a_2=a_3$ by (\ref{dPcondition}),
\item or $a_0+a_3\divby a_2$, i.e. $a_3=2a_2-a_0$ by (\ref{dPcondition}),
\item or $a_0+a_3=a_e+\mu a_2$, $3a_3=a_f+\nu a_2$ for some $e\neq f$, i.e. $a_1=4a_0/3$, $a_3=a_2+a_0/3$.
\end{itemize}

In the last case, $d_1=a_2+4a_0/3$, $d_2=3a_2+a_0$, $a_2<5a_0/3$. By Theorem~$WF$, $9a_2\divby a_0$, i.e. $a_2=\lambda a_0/9$, $\lambda <15$, $a_3=(\lambda +3)a_0/9$, and so $a_4=(\lambda +3)a_0/6 < 54$ by Lemma~$2$.\\

If $a_3=2a_2-a_0$, then $d_1=2a_2$, $d_2=6a_2-3a_0$, $a_2<a_1/2+3a_0/4$. By Theorem~$QS$~$(1)$, $i=1$, 
\begin{itemize}
\item either $2a_2\divby a_1$, i.e. $a_1=a_2$ by (\ref{dPcondition}),
\item or $6a_2-3a_0\divby a_1$, i.e. $a_2=2a_1/3+a_0/2$ by (\ref{dPcondition}).
\end{itemize}

The first case contradicts to (\ref{AllNonequal}). Hence $a_2=2a_1/3+a_0/2$, and so $d_1=4a_1/3+a_0$, $d_2=4a_1$, $a_3=4a_1/3$, $a_4=2a_1$. By Theorem~$WF$, $a_0\divby \frac{a_1}{3}$. Then $a_4\leq 12$ by Lemma~$2$.\\

If $a_4=a_3+a_2/2$, then $d_1=a_0+a_3$, $d_2=2a_3+a_2$, $a_3<a_1+a_2/2$. By Theorem~$QS$~$(1)$, $i=2$, 
\begin{itemize}
\item either $2a_3\divby a_2$, i.e. $a_2=a_3$ by (\ref{dPcondition}),
\item or $a_0+a_3\divby a_2$, i.e. $a_3=2a_2-a_0$ by (\ref{dPcondition}),
\item or $a_0+a_3=a_e+\mu a_2$, $2a_3=a_f+\nu a_2$ for some $e\neq f$, i.e. $a_1=3a_0/2$, $a_3=a_2+a_0/2$.
\end{itemize}

In the second case $a_4<320$, while in the third case $a_4<51$ by Theorem~$QS$, Theorem~$WF$ and Lemma~$2$.\\

If $a_4=a_3+a_1/2$, then $d_1=a_0+a_3$, $d_2=2a_3+a_1$, $a_3<a_2+a_1/2$. By Theorem~$QS$~$(1)$, $i=2$, 
\begin{itemize}
\item either $2a_3+a_1\divby a_2$, i.e. $a_3=3a_2/2-a_1/2$ by (\ref{dPcondition}),
\item or $a_0+a_3\divby a_2$, i.e. $a_3=2a_2-a_0$ by (\ref{dPcondition}),
\item or $a_0+a_3=a_e+\mu a_2$, $2a_3+a_1=a_f+\nu a_2$ for some $e\neq f$, i.e. ($a_2=3(a_1-a_0)$, $a_3=4(a_1-a_0)$) or ($a_2=2a_1-a_0$, $a_3=3a_1-2a_0$).
\end{itemize}

In the second case $a_4<5$, while in the third item $a_4<35$ by Theorem~$QS$, Theorem~$WF$ and Lemma~$2$.\\

If $a_3=3a_2/2-a_1/2$, then $d_1=3a_2/2-a_1/2+a_0$, $d_2=3a_2$, $a_2<2a_1$. By Theorem~$QS$~$(1)$, $i=1$, 
\begin{itemize}
\item either $3a_2\divby a_1$, i.e. $a_2=5a_1/3$ or $a_2=4a_1/3$ by (\ref{dPcondition}),
\item or $3a_2+2a_0\divby a_1$, i.e. $a_2=\lambda a_1/3-2a_0/3$, $\lambda \in \{ 4,5,6,7 \}$ by (\ref{dPcondition}),
\item or $3a_2/2-a_1/2+a_0=a_e+\mu a_1$, $3a_2=a_f+\nu a_1$ for some $e\neq f$.
\end{itemize}

In the first item $a_4<16$, while in the third item $a_4<64$ by Lemma~$2$ by the same analysis as in the previous cases.\\

If $a_2=7a_1/3-2a_0/3$, then $d_1=3a_1$, $d_2=7a_1-2a_0$, $a_1<2a_0$. Then $a_4<168$ by Theorem~$QS$~$(1)$, $i=0$, and Lemma~$2$.\\

By Theorem~$QS$~$(2)$, $(i,j)=(1,2)$, $d_1= 3a_2/2-a_1/2+a_0=a_e+\alpha a_1+\beta a_2$, $\alpha , \beta \geq 0$, for some $e$. Together with (\ref{dPcondition}), this implies that
\begin{itemize}
\item either $d_1=a_1+a_2$, i.e. $a_2=3a_1-2a_0$,
\item or $d_1=a_0+2a_1$, i.e. $a_2=5a_1/3$,
\item or $d_1=2a_1$, i.e. $a_2=5a_1/3-2a_0/3$,
\item or $d_1=3a_1$, i.e. $a_2=7a_1/3-2a_0/3$.
\end{itemize}

If $a_2=3a_1-2a_0=\lambda a_1/3-2a_0/3$, then $a_4<18$ by Lemma~$2$. Hence it remains to consider the following case:
$$
a_2=\frac{5a_1}{3}-\frac{2a_0}{3},\; a_3=2a_1-a_0,\; a_4=\frac{5a_1}{2}-a_0,\; d_1=2a_1,\; d_2=5a_1-2a_0. 
$$

By Theorem~$QS$~$(1)$, $i=0$, either $4a_1\divby a_0$ or $10a_1\divby a_0$. Since $5 \nmid a_0$ by Theorem~$WF$, $4a_1\divby a_0$. Hence
$$
a_1=\frac{\lambda a_0}{4},\; a_2=\left( \frac{5\lambda}{12}-\frac{2}{3} \right) a_0,\; a_3=\left( \frac{\lambda}{2}-1 \right) a_0,\; a_4=\left( \frac{5\lambda}{8}-1 \right) a_0.
$$

By Theorem~$WF$,
$$
a_0=1,\; a_1=2t,\; a_2=\frac{2}{3}(5t-1),\; a_3=4t-1,\; a_4=5t-1,\; d_1=4t,\; d_2=10t-2.
$$

This violates the condition of Theorem~$QS$~$(2)$, $(i,j)=(2,4)$.\\

{\it Case $(4.2)$:} $d_1=a_1+a_3$, $d_2=2a_4$, $a_4<a_0+a_2$ by (\ref{dPcondition}). By Theorem~$QS$~$(1)$, $i=3$, 
\begin{itemize}
\item either $2a_4\divby a_3$, i.e. $a_4=3a_3/2$ by (\ref{dPcondition}),
\item or $2a_4-a_0\divby a_3$, i.e. $a_4=a_3+a_0/2$ by (\ref{dPcondition}),
\item or $2a_4-a_2\divby a_3$, i.e. $a_4=a_3+a_2/2$ by (\ref{dPcondition}).
\end{itemize}

In the first case $a_4<23$, while in the third case $a_4<39$ by Theorem~$QS$ and Lemma~$2$.\\

If $a_4=a_3+a_0/2$, then $d_1=a_1+a_3$, $d_2=2a_3+a_0$, $a_3<a_2+a_0/2$. By Theorem~$QS$~$(1)$, $i=2$, 
\begin{itemize}
\item either $a_1+a_3\divby a_2$, i.e. $a_3=2a_2-a_1$ by (\ref{dPcondition}),
\item or $a_0+2a_3\divby a_2$, i.e. $a_3=3a_2/2-a_0/2$ by (\ref{dPcondition}).
\end{itemize}

In the first case $a_4<120$, while in the second case $a_4<36$ by Theorem~$QS$, Theorem~$WF$ and Lemma~$2$.\\

{\it Cases $(4.3)-(4.9)$:} By the same analysis, $a_4<75$. {\it QED}\\

\vspace{3ex}

{\bf Lemma~$6$.} {\it The} Main Theorem {\it holds in cases $(03)$, $(13)$ and $(23)$ of Lemma~$1$.} \\

{\it Proof:} We assume (\ref{AllNonequal}) and consider each case separately.\\

{\it Case $(03)$:} $d_1=a_0+a_4$, $d_2=a_3+a_4$, $a_4<a_1+a_2$ by (\ref{dPcondition}). By Theorem~$QS$~$(1)$, $i=3$, 
\begin{itemize}
\item either $a_4\divby a_3$, i.e. $a_4=a_3$ by (\ref{dPcondition}),
\item or $a_0+a_4\divby a_3$, i.e. $a_4=2a_3-a_0$ by (\ref{dPcondition}),
\item or $a_0+a_4-a_1\divby a_3$, i.e. $a_4=a_3+a_1-a_0$ by (\ref{dPcondition}),
\item or $a_0+a_4-a_2\divby a_3$, i.e. $a_4=a_3+a_2-a_0$ by (\ref{dPcondition}).
\end{itemize}

If $a_4=2a_3-a_0$, then $a_3=a_2+a_0/3$ by Theorem~$QS$~$(1)$, $i=2$. Then $d_1=2a_2+2a_0/3$, $d_2=3a_2$, $a_2<a_1+a_0/3$. By Theorem~$QS$~$(1)$, $i=1$, 
\begin{itemize}
\item either $d_1\divby a_1$, i.e. $a_2=3a_1/2-a_0/3$ by (\ref{dPcondition}),
\item or $d_1-a_0\divby a_1$, i.e. $a_2=a_1+a_0/6$ by (\ref{dPcondition}).
\end{itemize}

In the first case, $a_4<54$ by Theorem~$QS$ and Lemma~$2$.\\

In the second case, $a_4=2a_1$, $d_1=2a_1+a_0$, $d_2=3a_1+a_0/2$. By Theorem~$WF$, $a_1=a_0$, which contradicts to (\ref{AllNonequal}).\\

If $a_4=a_3+a_1-a_0$, then $d_1=a_1+a_3$, $d_2=2a_3+a_1-a_0$, $a_3<a_0+a_2$. By Theorem~$QS$~$(1)$, $i=2$, 
\begin{itemize}
\item either $a_1+a_3\divby a_2$, i.e. $a_3=2a_2-a_1$ by (\ref{dPcondition}),
\item or $2a_3+a_1-a_0\divby a_2$, i.e. $a_3=(3a_2+a_0-a_1)/2$ by (\ref{dPcondition}).
\end{itemize}

If $a_3=2a_2-a_1$, then $d_1=2a_2$, $d_2=4a_2-a_1-a_0$, $a_2<a_0+a_1$. By Theorem~$QS$~$(1)$, $i=1$, 
\begin{itemize}
\item either $2a_2\divby a_1$, i.e. $a_2=3a_1/2$ by (\ref{dPcondition}),
\item or $4a_2-a_0\divby a_1$, i.e. $a_2=\lambda a_1/4+a_0/4$, $\lambda \in \{ 4,5,6 \}$ by (\ref{dPcondition}),
\item or $4a_2-2a_0\divby a_1$, i.e. $a_2=\lambda a_1/4+a_0/2$, $\lambda \in \{ 3,4,5 \}$ by (\ref{dPcondition}),
\item or $3a_2-a_0\divby a_1$, i.e. $a_2=(\lambda a_1+a_0)/3$, $\lambda \in \{ 3,4 \}$ by (\ref{dPcondition}).
\end{itemize}

Theorem~$QS$~$(1)$, $i=0$, and Lemma~$2$ give bound $a_4<195$ in all cases except for the following two:
\begin{itemize}
\item $a_2=a_1+a_0/4$,
\item $a_2=a_1+a_0/3$.
\end{itemize}

In these cases, by Theorem~$QS$ and Theorem~$WF$,
\begin{itemize}
\item either $a_0=8$, $a_1=4t+1$, $a_2=4t+3$, $a_3=4t+5$, $a_4=8t-2$, $d_1=8t+6$, $d_2=12t+3$,
\item or $a_0=9$, $a_1=t$, $a_2=t+3$, $a_3=t+6$, $a_4=2t-3$, $d_1=2t+6$, $d_2=3t+3$, where $t\equiv -1 \pmod 3$.
\end{itemize}

These solutions appear in Table~\ref{table:1} (No.~$43$ and No.~$45$).\\

If $a_3=(3a_2+a_0-a_1)/2$, then $d_1=(3a_2+a_0+a_1)/2$, $d_2=3a_2$, $a_2<a_0+a_1$. By Theorem~$QS$~$(1)$, $i=1$, 
\begin{itemize}
\item either $3a_2\divby a_1$, i.e. $a_2=4a_1/3$ or $a_2=5a_1/3$ by (\ref{dPcondition}),
\item or $3a_2+a_0\divby a_1$, i.e. $a_2=\lambda a_1/3-a_0/3$, $\lambda \in \{ 4,5,6 \}$ by (\ref{dPcondition}),
\item or $3a_2-a_0\divby a_1$, i.e. $a_2=\lambda a_1/3+a_0/3$, $\lambda \in \{ 3,4 \}$ by (\ref{dPcondition}),
\item or $2a_2\divby a_1$, i.e. $a_2=3a_1/2$ by (\ref{dPcondition}).
\end{itemize}

In each case, $a_4<220$ by Theorem~$WF$, Theorem~$QS$ and Lemma~$2$.\\

If $a_4=a_3+a_2-a_0$, then $d_1=a_2+a_3$, $d_2=2a_3+a_2-a_0$, $a_3<a_0+a_1$. By Theorem~$QS$~$(1)$, $i=2$, $2a_3-a_0\divby a_2$, i.e. $a_3=a_2+a_0/2$ by (\ref{dPcondition}). Then $d_1=2a_2+a_0/2$, $d_2=3a_2$, $a_2<a_1+a_0/2$. By Theorem~$QS$~$(1)$, $i=1$, 
\begin{itemize}
\item either $3a_2\divby a_1$, i.e. $a_2=4a_1/3$ by (\ref{dPcondition}),
\item or $4a_2+a_0\divby a_1$, i.e. $a_2=\lambda a_1/4-a_0/4$, $\lambda \in \{ 5,6 \}$ by (\ref{dPcondition}),
\item or $2a_2-a_0/2\divby a_1$, i.e. $a_2=a_1+a_0/4$ by (\ref{dPcondition}).
\end{itemize}

In the first two items, $a_4<286$ by Theorem~$QS$ and Lemma~$2$. In the third item, $a_4=2a_1$, $d_1=2a_1+a_0$, $d_2=3a_1+3a_0/4$. By Theorem~$WF$, $3a_0\divby a_1$, i.e. $a_1=3a_0$ or $a_1=3a_0/2$. In each case, $a_4<25$ by Lemma~$2$.\\

{\it Case $(13)$:} $d_1=a_1+a_4$, $d_2=a_3+a_4$, $a_4<a_0+a_2$ by (\ref{dPcondition}). By Theorem~$QS$~$(1)$, $i=3$, 
\begin{itemize}
\item either $a_1+a_4\divby a_3$, i.e. $a_4=2a_3-a_1$ by (\ref{dPcondition}),
\item or $a_1+a_4-a_2\divby a_3$, i.e. $a_4=a_3+a_2-a_1$ by (\ref{dPcondition}).
\end{itemize}

If $a_4=2a_3-a_1$, then $a_3=a_2+a_1/3$ by Theorem~$QS$~$(1)$, $i=2$. Then $d_1=2a_2+2a_1/3$, $d_2=3a_2$, $a_2<a_0+a_1/3$. By Theorem~$QS$~$(1)$, $i=1$, $6a_2\divby a_1$, i.e. $a_2=7a_1/6$, $a_3=3a_1/2$, $a_4=2a_1\leq 12$ by Lemma~$2$.\\

If $a_4=a_3+a_2-a_1$, then $d_1=a_2+a_3$, $d_2=2a_3+a_2-a_1$, $a_3<a_0+a_1$. By Theorem~$QS$~$(1)$, $i=2$, 
\begin{itemize}
\item either $2a_3-a_1\divby a_2$, i.e. $a_3=a_2+a_1/2$ by (\ref{dPcondition}),
\item or $2a_3-a_1-a_0\divby a_2$, i.e. $a_3=(a_0+a_1+a_2)/2$ by (\ref{dPcondition}),
\item or $2a_3-2a_1\divby a_2$, i.e. $a_3=a_1+a_2/2$ by (\ref{dPcondition}).
\end{itemize}

In each case, $a_4<288$ by Theorem~$WF$, Theorem~$QS$ and Lemma~$2$.\\

{\it Case $(23)$:} $d_1=a_2+a_4$, $d_2=a_3+a_4$, $a_4<a_0+a_1$ by (\ref{dPcondition}). By Theorem~$QS$~$(1)$, $i=3$, $a_2+a_4\divby a_3$, i.e. $a_4=2a_3-a_2$. Then $d_1=2a_3$, $d_2=3a_3-a_2$, $a_3<(a_0+a_1+a_2)/2$. By Theorem~$QS$~$(1)$, $i=2$,
\begin{itemize}
\item either $3a_3\divby a_2$, i.e. $a_3=4a_2/3$ by (\ref{dPcondition}),
\item or $3a_3-a_0\divby a_2$, i.e. $a_3=a_2+a_0/3$ by (\ref{dPcondition}),
\item or $3a_3-a_1\divby a_2$, i.e. $a_3=a_2+a_1/3$ by (\ref{dPcondition}).
\end{itemize}

In each case, $a_4<42$ by Theorem~$WF$, Theorem~$QS$ and Lemma~$2$. {\it QED}\\

\vspace{3ex}

{\bf Lemma~$7$.} {\it Assume given integers $1\leq a_0 <a_1<a_2<a_3<a_4<2a_3$. Then $\forall \; i\neq j \in \{ 0,1,2 \}$,
\begin{align*}
a_i +a_4\in (a_3,a_4) \quad &\mbox{if and only if} \quad a_4=2a_3-a_i,\\
a_i+a_4-a_j\in (a_3,a_4) \quad &\mbox{if and only if} \quad a_4=\begin{cases}
2a_3+a_j-a_i\; &\text{if $i>j$},\\
a_3+a_j-a_i\; &\text{if $i<j$}.
\end{cases}
\end{align*}
Here $(a_3,a_4)=\{ \alpha a_3+\beta a_4 \;\mid \; \alpha,\beta \in \mathbb Z,\; \alpha, \beta \geq 0 \}$.}\\

{\it Proof:} If $a_i+a_4=\alpha a_3+\beta a_4$, then $\beta =0$ and then $\alpha a_3=a_i+a_4<3a_3$. Hence $\alpha =2$.\\

If $a_i+a_4=a_j+\alpha a_3+\beta a_4$, then $\beta =0$. Then $\alpha a_3=a_i-a_j+a_4<3a_3$ if $i>j$, and $\alpha a_3<a_4<2a_3$ if $i<j$. Hence $\alpha =2$ if $i>j$, and $\alpha =1$ if $i<j$. {\it QED}\\

\vspace{3ex}

{\bf Corollary~$2$.} {\it The condition of Theorem~$QS$~$(2)$, $(i,j)=(3,4)$, can be rewritten as follows in cases $(01)$, $(02)$ and $(12)$ of Lemma~$1$, assuming $a_i\neq a_j$ for $i\neq j$ (and given the assumptions of Lemma~$1$):
\begin{itemize}
\item (01): \begin{itemize}
\item either $a_4=2a_3-a_0$,
\item or $a_4=2a_3-a_1$,
\item or $a_2=2a_1-a_0$, $a_4=a_3+a_1-a_0$,
\item or $a_3=a_2+a_1-2a_0$, $a_4=2a_2+a_1-3a_0$;
\end{itemize}
\item (02): \begin{itemize}
\item either $a_4=2a_3-a_0$,
\item or $a_4=2a_3-a_2$,
\item or $a_3=a_2-a_0$, $a_4=a_2+a_1-2a_0$,
\item or $a_3=a_2+a_1-2a_0$, $a_4=a_2+2a_1-3a_0$,
\item or $a_3=2a_2-a_1-a_0$, $a_4=3a_2-a_1-2a_0$;
\end{itemize}
\item (12): \begin{itemize}
\item either $a_4=2a_3-a_1$,
\item or $a_4=2a_3-a_2$,
\item or $a_2=2a_1-a_0$, $a_4=2a_3+a_0-a_1$,
\item or $a_3=2a_2-a_1-a_0$, $a_4=3a_2-2a_1-a_0$.
\end{itemize}
\end{itemize}}

{\it Proof:} Let us consider Case $(01)$. The other two cases are similar. The condition 
$$
d_1-a_e, \; d_2-a_f\in (a_3,a_4),\;\; \#\{e\}=\#\{f\}=2,\; \{e,f\}=\{0,1,2\}
$$ 
of Theorem~$QS$ is the same as:
\begin{itemize} 
\item either $d_1-a_2, \; d_2-a_2\in (a_3,a_4)$,
\item or $d_1-a_1, \; d_2-a_2\in (a_3,a_4)$,
\item or $d_1-a_2, \; d_2-a_0\in (a_3,a_4)$.
\end{itemize}

The first case can not occur, because $a_0\neq a_1$. Hence by Lemma~$7$,
\begin{itemize} 
\item either $a_4=a_3+a_1-a_0=a_3+a_2-a_1$,
\item or $a_4=a_3+a_2-a_0=2a_3+a_0-a_1$.
\end{itemize}
{\it QED}\\

\section{Proof of the {\bf Main Theorem}}

By Lemma~$5$ and Lemma~$6$, it remains to consider cases $(01)$, $(02)$ and $(12)$ of Lemma~$1$. Each of these cases splits into several subcases according to Corollary~$2$. Lemma~$4$ allows us to assume (\ref{AllNonequal}).

\subsection{Case $(01)$}

\subsubsection{Subcase $a_4=2a_3-a_0$}

In this subcase, $d_1=2a_3$, $d_2=2a_3+a_1-a_0$, $a_3<a_0+a_2$. By Theorem~$QS$~$(1)$, $i=2$, 
\begin{itemize}
\item either $2a_3\divby a_2$, i.e. $a_3=3a_2/2$ by (\ref{dPcondition}),
\item or $2a_3+a_1-a_0\divby a_2$, i.e. $a_3=(3a_2+a_0-a_1)/2$ by (\ref{dPcondition}),
\item or $2a_3-a_0\divby a_2$, i.e. $a_3=a_2+a_0/2$ by (\ref{dPcondition}).
\end{itemize}

In each case, $a_4<198$ by Theorem~$QS$, Theorem~$WF$ and Lemma~$2$.\\

\subsubsection{Subcase $a_4=2a_3-a_1$}

In this subcase, $d_1=2a_3+a_0-a_1$, $d_2=2a_3$, $a_3<a_1+a_2$. By Theorem~$QS$~$(1)$, $i=2$, 
\begin{itemize}
\item either $2a_3\divby a_2$, i.e. $a_3=3a_2/2$ by (\ref{dPcondition}),
\item or $2a_3+a_0-a_1\divby a_2$, i.e. $a_3=(3a_2+a_1-a_0)/2$ or $a_3=a_2+(a_1-a_0)/2$ by (\ref{dPcondition}),
\item or $a_3+a_0-a_1\divby a_2$, i.e. ($a_1=2a_0$, $a_3=a_2+a_0$) or ($a_1=3a_0/2$, $a_3=a_2+a_0/2$) or ($a_2=2a_1-3a_0$, $a_3=3a_1-4a_0$) by (\ref{dPcondition}),
\item or $2a_3+a_0-2a_1\divby a_2$, i.e. $a_2=2(a_1-a_0)$, $a_3=3a_1-5a_0/2$ by (\ref{dPcondition}),
\item or $2a_3-a_1\divby a_2$, i.e. $a_3=a_2+a_1/2$ by (\ref{dPcondition}).
\end{itemize}

If $a_3=a_2+a_1/2$, then $a_4=2a_2$, $d_1=2a_2+a_0$, $d_2=2a_2+a_1$. By Theorem~$WF$, $a_1=a_2$, which contradicts to (\ref{AllNonequal}).\\

If $a_2=2(a_1-a_0)$, $a_3=3a_1-5a_0/2$, then $a_4=5(a_1-a_0)$, $d_1=a_0+5(a_1-a_0)$, $d_2=a_0+6(a_1-a_0)$. Hence $a_0\divby (a_1-a_0)$ by Theorem~$WF$, i.e. $a_1=\frac{\lambda +1}{\lambda} a_0$. Then $a_2=\frac{2}{\lambda} a_0$, and so $a_4=\frac{5}{\lambda}a_0\leq 5$ by Lemma~$2$.\\

If $a_1=2a_0$, $a_3=a_2+a_0$, then $a_4=2a_2$, $d_1=2a_2+a_0$, $d_2=2a_2+2a_0$. By Theorem~$WF$, $2a_0\divby a_2$, i.e. $a_2=2a_0=a_1$, which contradicts to (\ref{AllNonequal}).\\

If $a_2=2a_1-3a_0$, $a_3=3a_1-4a_0$, then $d_1=5a_1-7a_0$, $d_2=6a_1-8a_0$. By Theorem~$QS$~$(1)$, $i=1$, either $3a_0\divby a_1$ or $7a_0\divby a_1$ or $8a_0\divby a_1$, which implies that $a_4<40$ by Lemma~$2$.\\

If $a_1=3a_0/2$, $a_3=a_2+a_0/2$, then $d_1=2a_2+a_0/2$, $d_2=2a_2+a_0$. By Theorem~$WF$, $2a_2\divby a_0/2$, i.e. 
$$
a_2=\frac{\lambda a_0}{4},\; a_3=\frac{\lambda +2}{4}a_0,\; a_4=\frac{\lambda -1}{2}a_0,\; d_1=\frac{\lambda +1}{2}a_0,\; d_2=\frac{\lambda +2}{2}a_0.
$$

Hence \begin{itemize}
\item either $a_0=4$, $a_1=6$, $a_2=t$, $a_3=t+2$, $a_4=2t-2$, $d_1=2t+2$, $d_2=2t+4$, where $t\equiv -1\pmod{3}$ is odd,
\item or $a_0=2$, $a_1=3$, $a_2=t$, $a_3=t+1$, $a_4=2t-1$, $d_1=2t+1$, $d_2=2t+2$, where $t\equiv 1\pmod{3}$.
\end{itemize}

These solutions appear in Table~\ref{table:1} (No.~$23$ and No.~$31$).\\

If $a_3=3a_2/2$, then $d_1=3a_2+a_0-a_1$, $d_2=3a_2$, $a_2<2a_1$. By Theorem~$QS$~$(2)$, $(i,j)=(2,3)$, 
\begin{itemize}
\item either $a_1-a_0\divby a_2/2$, i.e. $a_2=2(a_1-a_0)$,
\item or $a_0\divby a_2/2$, i.e. $a_2=2a_0$,
\item or $2a_1-a_0\divby a_2/2$, i.e. $a_2=2a_1-a_0$ or $a_2=(4a_1-2a_0)/3$.
\end{itemize}

If $a_2=2(a_1-a_0)$, then $d_1=5(a_1-a_0)$, $d_2=6(a_1-a_0)$. By Theorem~$QS$~$(1)$, $i=1$, either $5a_0\divby a_1$ or $6a_0\divby a_1$. Then $a_4<30$ by Lemma~$2$.\\

If $a_2=2a_0$, then $d_1=7a_0-a_1$, $d_2=6a_0$. By Theorem~$QS$~$(1)$, $i=1$, $N\cdot a_0\divby a_1$, where $N\in \{ 4,5,6,7 \}$. Then $a_4<35$ by Lemma~$2$.\\

If $a_2=2a_1-a_0$, then $d_1=5a_1-2a_0$, $d_2=6a_1-3a_0$. By Theorem~$QS$~$(1)$, $i=1$, either $2a_0\divby a_1$ or $3a_0\divby a_1$. Then $a_4<15$ by Lemma~$2$.\\

If $a_2=(4a_1-2a_0)/3$, then $d_1=3a_1-a_0$, $d_2=4a_1-2a_0$. By Theorem~$QS$~$(1)$, $i=0$, either $3a_1\divby a_0$ or $8a_1\divby a_0$. By Theorem~$WF$, $a_0=1$, i.e.
$$
a_0=1,\; a_1=2t+1,\; a_2=\frac{8t+2}{3},\; a_3=4t+1,\; a_4=6t+1,\; d_1=6t+2,\; d_2=8t+2.
$$

This solution appears in Table~\ref{table:1} (No.~$20$).\\

If $a_3=(3a_2+a_1-a_0)/2$, then $d_1=3a_2$, $d_2=3a_2+a_1-a_0$, $a_2<a_0+a_1$. By Theorem~$QS$~$(1)$, $i=1$, 
\begin{itemize}
\item either $3a_2\divby a_1$, i.e. $a_2=4a_1/3$ or $a_2=5a_1/3$ by (\ref{dPcondition}),
\item or $3a_2-a_0\divby a_1$, i.e. $a_2=(\lambda a_1+a_0)/3$, $\lambda \in \{ 3,4 \}$ by (\ref{dPcondition}),
\item or $d_1-a_3\divby a_1$, i.e. $a_2=(\lambda a_1-a_0)/3$, $\lambda \in \{ 4,5,6 \}$ by (\ref{dPcondition}),
\item or $2a_2\divby a_1$, i.e. $a_2=3a_1/2$ by (\ref{dPcondition}).
\end{itemize}

In each case, $a_4<308$ by Theorem~$WF$, Theorem~$QS$ and Lemma~$2$.\\

If $a_3=a_2+(a_1-a_0)/2$, then $d_1=2a_2$, $d_2=2a_2+a_1-a_0$. By Theorem~$QS$~$(1)$, $i=1$, 
\begin{itemize}
\item either $2a_2\divby a_1$, i.e. $a_2=\lambda a_1/2$, $a_3=(\lambda +1)a_1/2-a_0/2$, $a_4=\lambda a_1-a_0$, $d_1=\lambda a_1$, $d_2=(\lambda +1)a_1-a_0$,
\item or $2a_2-a_0\divby a_1$, i.e. $a_2=(\lambda a_1+a_0)/2$, $a_3=(\lambda +1)a_1/2$, $a_4=\lambda a_1$, $d_1=\lambda a_1+a_0$, $d_2=(\lambda +1)a_1$,
\item or $2a_2-a_3\divby a_1$, i.e. $a_2=\lambda a_1+(a_1-a_0)/2$, $a_3=(\lambda +1)a_1-a_0$, $a_4=(2\lambda +1) a_1-2a_0$, $d_1=(2\lambda +1) a_1-a_0$, $d_2=(2\lambda +2)a_1-2a_0$.
\end{itemize}

In the second case, $2a_0\divby a_1$ by Theorem~$WF$. Hence 
$$
a_1=2a_0,\; a_2=\left( \lambda +\frac{1}{2} \right) a_0,\; a_3=(\lambda +1) a_0,\; a_4=2\lambda a_0.
$$

This violates conditions of Theorem~$WF$.\\

In the first case, $\gcd(a_0,a_1) =1$ if $\lambda $ is even, and $\gcd(a_0,a_1) =2$ if $\lambda $ is odd. In the third case, $\gcd(a_0,a_1) =1$ by Theorem~$WF$. In these two cases, by Theorem~$QS$~$(1)$, $i=0$, \begin{itemize}
\item either $\lambda a_1 \divby a_0$ or $(\lambda +1) a_1 \divby a_0$ or ($(\lambda /2+1) \divby a_0$ and $\lambda$ is even) in the first case;
\item and either $(2\lambda +1) \divby a_0$ or $(\lambda +1) \divby a_0$ or $(2\lambda +3) \divby a_0$ in the third case.
\end{itemize}

If $\lambda a_1 \divby a_0$ in the first case, then $d_2-a_1=(\lambda a_1/a_0-1)a_0\in (a_0,a_2,a_4)$. Hence by Theorem~$QS$~$(3)$, $(k,l,m)=(0,2,4)$, either $d_2=2a_3\in (a_0,a_2,a_4)$ or $a_3=d_2-a_3\in (a_0,a_2,a_4)$, i.e. $2a_3\in (a_0,a_2,a_4)$. Hence $2a_3\in (a_0,a_2)$. This implies that either $(\lambda +1) a_1 \divby a_0$ or ($(\lambda /2+1) \divby a_0$ and $\lambda$ is even).\\

If $(2\lambda +1) \divby a_0$ in the third case, then $d_2-a_1=a_4\in (a_0,a_2,a_4)$. Hence by Theorem~$QS$~$(3)$, $(k,l,m)=(0,2,4)$, $d_2=2a_3\in (a_0,a_2,a_4)$. Hence $2a_3\in (a_0,a_2)$. This implies that either $(\lambda +1) \divby a_0$ or $(2\lambda +3) \divby a_0$.\\

Hence \begin{itemize}
\item either $a_0=2b_0$, $a_1=2b_1$, $a_2=(\nu b_0-1)b_1$, $a_3=(\nu b_1-1)b_0$, $a_4=2(\nu b_0b_1-b_0-b_1)$, $d_1=2b_1(\nu b_0-1)$, $d_2=2b_0(\nu b_1-1)$, where $\gcd(b_0, b_1)=1$, $b_0$, $b_1$ are odd, $\nu $ is even;
\item or $a_2=(\nu a_0-1)a_1/2$, $a_3=(\nu a_1-1)a_0/2$, $a_4=\nu a_0a_1-a_0-a_1$, $d_1=a_1(\nu a_0-1)$, $d_2=a_0(\nu a_1-1)$, where $\gcd(a_0,a_1)=1$, $\nu$, $a_0$, $a_1$ are odd;
\item or $a_2=(\nu a_0-1)a_1$, $a_3=\nu a_0a_1 -(a_0+a_1)/2$, $a_4=2\nu a_0a_1-a_0-2a_1$, $d_1=2a_1(\nu a_0-1)$, $d_2=2\nu a_0a_1-a_0-a_1$, where $\gcd(a_0,a_1)=1$, $a_0$, $a_1$ are odd;
\item or $a_2=(\nu a_0-1)a_1+(a_1-a_0)/2$, $a_3=a_0(\nu a_1 -1)$, $a_4=2\nu a_0a_1-2a_0-a_1$, $d_1=2\nu a_0a_1-a_0-a_1$, $d_2=2a_0(\nu a_1-1)$, where $\gcd(a_0,a_1)=1$, $a_0$, $a_1$ are odd;
\item or $a_2=(\nu a_1-1)a_0/2-a_1$, $a_3=(\nu a_0-1)a_1/2-a_0$, $a_4=\nu a_0a_1-2(a_0+a_1)$, $d_1=\nu a_0a_1-a_0-2a_1$, $d_2=\nu a_0a_1-2a_0-a_1$, where $\gcd(a_0,a_1)=1$, $\nu$, $a_0$, $a_1$ are odd.
\end{itemize}

These solutions appear in Table~\ref{table:1} (Nos.~$10$-$14$).\\

\subsubsection{Subcase $a_2=2a_1-a_0$, $a_4=a_3+a_1-a_0$}

In this subcase, $d_1=a_3+a_1$, $d_2=a_3+2a_1-a_0=a_2+a_3$. By Theorem~$QS$~$(1)$, $i=2$, 
\begin{itemize}
\item either $a_3\divby a_2$, i.e. $a_3=(\mu-1)(2a_1-a_0)$,
\item or $a_3+a_1\divby a_2$, i.e. $a_3=\mu (2a_1-a_0)-a_1$,
\item or $a_1+a_3-a_0\divby a_2$, i.e. $a_4\divby a_2$.
\end{itemize}

If $a_4\divby a_2$, then by Theorem~$WF$, either $d_1=a_0+a_4\divby a_2$ or $d_2=a_1+a_4\divby a_2$, i.e. $a_1=a_2$. This contradicts to (\ref{AllNonequal}).\\

By Theorem~$QS$~$(1)$, $i=0$, 
\begin{itemize}
\item either $a_3+a_1\divby a_0$, i.e. $a_3=\lambda a_0-a_1$,
\item or $a_3+2a_1\divby a_0$, i.e. $a_3=\lambda a_0-2a_1$,
\item or $d_2-a_2\divby a_0$, i.e. $a_3\divby a_0$,
\item or $d_2-a_3\divby a_0$, i.e. $2a_1\divby a_0$, and hence $a_1\divby a_0$.
\end{itemize}

Note that $a_0$ is odd and $\gcd(a_0,a_1)=1$ by Theorem~$WF$. Moreover, $a_3\divby a_0$ if and only if $a_1\divby a_0$ if and only if $a_0=1$ by Theorem~$WF$.\\

If $a_3=\lambda a_0-2a_1=(\mu -1)(2a_1-a_0)$, then $\mu=\nu a_0$,
\begin{align*}
a_2=2a_1-a_0,\; a_3=(\nu a_0-1)(2a_1-a_0), \; a_4=\nu a_0 (2a_1-a_0)-a_1,\; d_1=(\nu a_0-1)(2a_1-a_0)+a_1,\\
d_2=\nu a_0 (2a_1-a_0),\; \mbox{where}\; \gcd(\nu a_0-1,a_1-a_0)=1.
\end{align*}

If $a_3=\lambda a_0-2a_1=\mu (2a_1-a_0)-a_1$, then $\mu=(\nu a_0-1)/2$, $\nu$ is odd,
\begin{align*}
a_2=2a_1-a_0,\; a_3=\left( \frac{\nu a_0-1}{2}\right) (2a_1-a_0) -a_1, \; a_4=\left( \frac{\nu a_0-1}{2}\right) (2a_1-a_0) -a_0,\\
d_1=\left( \frac{\nu a_0-1}{2}\right) (2a_1-a_0),\; d_2=\left( \frac{\nu a_0+1}{2}\right) (2a_1-a_0) -a_1,\; \mbox{where}\; \gcd\left(\frac{\nu a_0-3}{2},a_1-a_0\right) =1.
\end{align*}

If $a_3=\lambda a_0-a_1=(\mu -1)(2a_1-a_0)$, then $\mu=(\nu a_0+1)/2$, $\nu$ is odd,
\begin{align*}
a_2=2a_1-a_0,\; a_3=\left( \frac{\nu a_0-1}{2}\right) (2a_1-a_0), \; a_4=\left( \frac{\nu a_0-1}{2}\right) (2a_1-a_0)+a_1 -a_0,\\
d_1=\left( \frac{\nu a_0-1}{2}\right) (2a_1-a_0)+a_1,\; d_2=\left( \frac{\nu a_0+1}{2}\right) (2a_1-a_0),\; \mbox{where}\; \gcd\left(\frac{\nu a_0-1}{2},a_1-a_0\right) =1.
\end{align*}

If $a_3=\lambda a_0-a_1=\mu (2a_1-a_0)-a_1$, then $\mu=\nu a_0$,
\begin{align*}
a_2=2a_1-a_0,\; a_3=\nu a_0 (2a_1-a_0)-a_1, \; a_4=\nu a_0 (2a_1-a_0)-a_0,\; d_1=\nu a_0(2a_1-a_0),\\
d_2=(\nu a_0+1) (2a_1-a_0)-a_1,\; \mbox{where}\; \gcd(\nu a_0-1,a_1-a_0)=1.
\end{align*}

In all these cases $a_0$ is odd and $\gcd(a_0,a_1)=1$.\\

These solutions appear in Table~\ref{table:1} (Nos.~$6$-$9$).\\

\subsubsection{Subcase $a_3=a_2+a_1-2a_0$, $a_4=2a_2+a_1-3a_0$}

In this subcase, $d_1=2a_2+a_1-2a_0$, $d_2=2a_2+2a_1-3a_0$. By Theorem~$QS$~$(1)$, $i=2$, 
\begin{itemize}
\item either $2a_1-4a_0\divby a_2$, i.e. $a_2=2a_1-4a_0$,
\item or $2a_1-3a_0\divby a_2$, i.e. $a_2=2a_1-3a_0$,
\item or $a_1-3a_0\divby a_2$, i.e. $a_1=3a_0$.
\end{itemize}

In each case, $a_4<55$ by Theorem~$WF$, Theorem~$QS$ and Lemma~$2$.\\

\subsection{Case $(02)$}

\subsubsection{Subcase $a_4=2a_3-a_0$}

In this subcase, $d_1=2a_3$, $d_2=2a_3+a_2-a_0$, $a_3<a_0+a_1$. By Theorem~$QS$~$(1)$, $i=2$, 
\begin{itemize}
\item either $2a_3\divby a_2$, i.e. $a_3=3a_2/2$ by (\ref{dPcondition}),
\item or $2a_3-a_0\divby a_2$, i.e. $a_3=a_2+a_0/2$ by (\ref{dPcondition}),
\item or $2a_3-a_1\divby a_2$, i.e. $a_3=a_2+a_1/2$ by (\ref{dPcondition}).
\end{itemize}

In each case, $a_4<90$ by Theorem~$QS$ and Lemma~$2$.\\

\subsubsection{Subcase $a_4=2a_3-a_2$}

In this subcase, $d_1=2a_3+a_0-a_2$, $d_2=2a_3$, $a_3<a_1+a_2$. By Theorem~$QS$~$(1)$, $i=2$, 
\begin{itemize}
\item either $2a_3\divby a_2$, i.e. $a_3=3a_2/2$ by (\ref{dPcondition}),
\item or $2a_3+a_0\divby a_2$, i.e. $a_3=(3a_2-a_0)/2$ or $a_3=2a_2-a_0/2$ by (\ref{dPcondition}),
\item or $a_3+a_0\divby a_2$, i.e. $a_3=2a_2-a_0$ by (\ref{dPcondition}),
\item or $d_1-a_1\divby a_2$, i.e. $a_3=a_2+(a_1-a_0)/2$ or $a_3=(3a_2+a_1-a_0)/2$ by (\ref{dPcondition}).
\end{itemize}

If $a_3=(3a_2+a_1-a_0)/2$, then $d_1=2a_2+a_1$, $d_2=3a_2+a_1-a_0$, $a_2<a_1+a_0$. By Theorem~$QS$~$(1)$, $i=1$, 
\begin{itemize}
\item either $2a_2\divby a_1$, i.e. $a_2=3a_1/2$ by (\ref{dPcondition}),
\item or $3a_2-a_0\divby a_1$, i.e. $a_2=a_1+a_0/3$ or $a_2=(4a_1+a_0)/3$ by (\ref{dPcondition}),
\item or $2a_2-a_0\divby a_1$, i.e. $a_2=a_1+a_0/2$ by (\ref{dPcondition}).
\end{itemize}

If $a_2=3a_1/2$ or $a_2=(4a_1+a_0)/3$, then $a_4<154$ by Lemma~$2$ and Theorem~$QS$.\\

If $a_2=a_1+a_0/3$, then $a_3=2a_1$, $a_4=3a_1-a_0/3$, $d_1=3a_1+2a_0/3$, $d_2=4a_1$. By Theorem~$QS$~$(2)$, $(i,j)=(1,3)$, $2a_0\divby a_1$, i.e. $a_1=2a_0$, $a_3=4a_0$, and so $a_4=17a_0/3\leq 17$ by Lemma~$2$.\\

If $a_2=a_1+a_0/2$, then $a_4=3a_1$, $d_1=3a_1+a_0$, $d_2=4a_1+a_0/2$. By Theorem~$WF$, $a_0\divby a_1$, i.e. $a_0=a_1$, which contradicts to (\ref{AllNonequal}).\\

If $a_3=a_2+(a_1-a_0)/2$, then $d_1=a_2+a_1$, $d_2=2a_2+a_1-a_0$. By Theorem~$QS$~$(1)$, $i=1$, 
\begin{itemize}
\item either $a_2\divby a_1$, i.e. $a_2=\lambda a_1$, $a_3=(\lambda +1/2)a_1-a_0/2$, $a_4=(\lambda +1)a_1-a_0$, $d_1=(\lambda +1)a_1$, $d_2=(2\lambda +1)a_1-a_0$,
\item or $2a_2-a_0\divby a_1$, i.e. $a_2=(\lambda a_1+a_0)/2$, $a_3=(\lambda +1)a_1/2$, $a_4=(\lambda +2)a_1/2-a_0/2$, $d_1=(\lambda +2)a_1/2+a_0/2$, $d_2=(\lambda +1)a_1$,
\item or $2a_2-2a_0\divby a_1$, i.e. $a_2=\lambda a_1/2+a_0$, $a_3=(\lambda +1)a_1/2+a_0/2$, $a_4=(\lambda +2)a_1/2$, $d_1=(\lambda +2)a_1/2+a_0$, $d_2=(\lambda +1)a_1+a_0$.
\end{itemize}

Note that in the first item, $\gcd(a_0,a_1)=1$, $a_0$, $a_1$ are odd. In the second item, $\gcd(a_0,a_1)\in \{ 1,2 \}$ if $\lambda $ is odd, and $\gcd(a_0,a_1)=2$ if $\lambda $ is even.\\

In the third item, $a_1=2a_0$ by Theorem~$QS$~$(2)$, $(i,j)=(1,4)$. Then $a_0$ should be even, which violates conditions of Theorem~$WF$.\\

By Theorem~$QS$~$(1)$, $i=0$, 
\begin{itemize}
\item either $\lambda +1 \divby a_0$ or $2\lambda +1 \divby a_0$ or $\lambda \divby a_0$ in the first item;
\item and either $(\lambda +1)a_1 \divby a_0$ or $\lambda a_1 \divby a_0$ or ($(\lambda +2)a_1 \divby a_0$ and $(\lambda +2)a_1/a_0$ is odd) in the second item.
\end{itemize}

In the first item, $\lambda \divby a_0$ implies that either $\lambda +1 \divby a_0$ or $2\lambda +1 \divby a_0$ by Theorem~$QS$~$(2)$, $(i,j)=(0,2)$. In the second item, $\lambda a_1 \divby a_0$ implies that either $(\lambda +1)a_1 \divby a_0$ or ($(\lambda +2)a_1 \divby a_0$ and $(\lambda +2)a_1/a_0$ is odd) by Theorem~$QS$~$(2)$, $(i,j)=(0,2)$.\\

Hence \begin{itemize}
\item either $a_2=(\nu a_0-1)a_1$, $a_3=\nu a_0a_1-(a_0+a_1)/2$, $a_4=(\nu a_1 -1)a_0$, $d_1=\nu a_0a_1$, $d_2=2\nu a_0a_1-a_0-a_1$, where $\gcd(a_0,a_1)=\gcd(\nu a_0-1,(a_1-a_0)/2)=1$, $a_0$, $a_1$ are odd;
\item or $a_2=(\nu a_0-1)a_1/2$, $a_3=(\nu a_1-1)a_0/2$, $a_4=(\nu a_0+1)a_1/2-a_0$, $d_1=(\nu a_0+1)a_1/2$, $d_2=(\nu a_1-1)a_0$, where $\gcd(a_0,a_1)=\gcd((\nu a_0-1)/2,(a_1-a_0)/2)=1$, $\nu$, $a_0$, $a_1$ are odd;
\item or $a_0=2b_0$, $a_1=2b_1$, $a_2=(\nu b_0-1)b_1+b_0$, $a_3=\nu b_0b_1$, $a_4=(\nu b_0+1)b_1-b_0$, $d_1=(\nu b_0+1)b_1+b_0$, $d_2=2\nu b_0b_1$, where $\gcd(b_0,b_1)=\gcd(\nu,b_1-b_0)=1$, $\nu$, $b_0$, $b_1$ are odd;
\item or $a_2=\nu a_0a_1+(a_0-a_1)/2$, $a_3=\nu a_0a_1$, $a_4=\nu a_0a_1+(a_1-a_0)/2$, $d_1=\nu a_0a_1+(a_0+a_1)/2$, $d_2=2\nu a_0a_1$, where $\gcd(a_0,a_1)=\gcd(\nu,(a_1-a_0)/2)=1$, $a_0$, $a_1$ are odd;
\item or $a_2=(\nu a_1+1)a_0/2-a_1$, $a_3=(\nu a_0-1)a_1/2$, $a_4=(\nu a_1-1)a_0/2$, $d_1=(\nu a_1+1)a_0/2$, $d_2=(\nu a_0-1)a_1$, where $\gcd(a_0,a_1)=\gcd((\nu a_0-1)/2,(a_1-a_0)/2)=1$, $\nu$, $a_0$, $a_1$ are odd.
\end{itemize}

These solutions appear in Table~\ref{table:1} (Nos.~$1$-$5$).\\

If $a_3=3a_2/2$, then $a_4=2a_2$, $d_1=2a_2+a_0$, $d_2=3a_2$. Hence $a_0\divby a_2/2$ by Theorem~$WF$, i.e. $a_2=2a_0$. Then $a_3=3a_0$, $a_4=4a_0\leq 4$ by Lemma~$2$.\\

If $a_3=(3a_2-a_0)/2$, then $d_1=2a_2$, $d_2=3a_2-a_0$, $a_2<2a_1+a_0$. By Theorem~$QS$~$(1)$, $i=1$, 
\begin{itemize}
\item either $2a_2\divby a_1$, i.e. $a_2=\lambda a_1/2$, $\lambda \in \{ 3,4,5 \}$ by (\ref{dPcondition}),
\item or $3a_2-a_0\divby a_1$, i.e. $a_2=(\lambda a_1+a_0)/3$, $3\leq \lambda \leq 7$ by (\ref{dPcondition}),
\item or $2a_2-a_0\divby a_1$, i.e. $a_2=(\lambda a_1+a_0)/2$, $\lambda \in \{ 2,3,4 \}$ by (\ref{dPcondition}),
\item or $3a_2-2a_0\divby a_1$ and $2a_2-a_3\divby a_1$, i.e. $a_2=2a_1-a_0$ and $a_1=5a_0/\mu$, $1\leq \mu \leq 4$ by (\ref{dPcondition}).
\end{itemize}

In the last item, $a_4<20$ by Lemma~$2$. In the second item, $4a_0\divby a_1$ by Theorem~$QS$~$(2)$, $(i,j)=(1,3)$. Then $a_4<56$ by Lemma~$2$. In the third item, $a_4=\lambda a_1$, $d_1=\lambda a_1+a_0$, $d_2=(3\lambda a_1+a_0)/2$. Hence $a_0=a_1$ by Theorem~$WF$, which violates (\ref{AllNonequal}).\\

If $a_2=5a_1/2$, then $a_4<450$ by Theorem~$QS$ and Lemma~$2$.\\

If $a_2=3a_1/2$, then $d_1=3a_1$, $d_2=9a_1/2-a_0$. By Theorem~$QS$~$(1)$, $i=0$, either $9a_1\divby a_0$ or $7a_1\divby a_0$. Since $3\nmid a_0$ by Theorem~$WF$, $7a_1\divby a_0$, i.e.
$$
a_1=\frac{\lambda a_0}{7},\; a_2=\frac{3\lambda a_0}{14},\; a_3=\left( \frac{9\lambda}{28}-\frac{1}{2} \right) a_0,\; a_4=\left( \frac{3\lambda}{7}-1 \right) a_0,\; d_1=\frac{3\lambda a_0}{7},\; d_2=\left( \frac{9\lambda}{14}-1 \right) a_0. 
$$

Hence \begin{itemize}
\item either $a_0=1$, $a_1=2t$, $a_2=3t$, $a_3=(9t-1)/2$, $a_4=6t-1$, $d_1=6t$, $d_2=9t-1$, where $t$ is odd,
\item or $a_0=7$, $a_1=2t$, $a_2=3t$, $a_3=(9t-7)/2$, $a_4=6t-7$, $d_1=6t$, $d_2=9t-7$, where $t$ is odd, $7\nmid t$.
\end{itemize}

These solutions appear in Table~\ref{table:1} (No.~$18$ and No.~$41$).\\

If $a_2=2a_1$, then $d_1=4a_1$, $d_2=6a_1-a_0$. By Theorem~$QS$~$(1)$, $i=0$, either $6a_1\divby a_0$ or $5a_1\divby a_0$ or $4a_1\divby a_0$. Since $\gcd(a_0,a_1)=1$ by Theorem~$WF$ and $a_0$ is even, $a_0\in \{ 2,4,6 \}$, i.e. \begin{itemize}
\item either $a_0=2$, $a_2=2a_1$, $a_3=3a_1-1$, $a_4=4a_1-2$, $d_1=4a_1$, $d_2=6a_1-2$,
\item or $a_0=4$, $a_2=2a_1$, $a_3=3a_1-2$, $a_4=4a_1-4$, $d_1=4a_1$, $d_2=6a_1-4$,
\item or $a_0=6$, $a_2=2a_1$, $a_3=3a_1-3$, $a_4=4a_1-6$, $d_1=4a_1$, $d_2=6a_1-6$.
\end{itemize}

In the first and the third items, either $a_1$ or $a_3$ is even, which violates conditions of Theorem~$WF$. Hence
$$
a_0=4,\; a_1=2t+1,\; a_2=4t+2,\; a_3=6t+1,\; a_4=8t,\; d_1=8t+4,\; d_2=12t+2.
$$

This solution appears in Table~\ref{table:1} (No.~$36$).\\

If $a_3=2a_2-a_0$, then $d_1=3a_2-a_0$, $d_2=4a_2-2a_0$, $a_2<a_0+a_1$. By Theorem~$QS$~$(1)$, $i=1$, 
\begin{itemize}
\item either $3a_2-a_0\divby a_1$, i.e. $a_2=a_1+a_0/3$ or $a_2=(4a_1+a_0)/3$ by (\ref{dPcondition}),
\item or $4a_2-2a_0\divby a_1$, i.e. $a_2=\lambda a_1/4+a_0/2$, $\lambda \in \{ 3,4,5 \}$ by (\ref{dPcondition}),
\item or $3a_2-2a_0\divby a_1$, i.e. $a_2=a_1+2a_0/3$ or $a_2=2(a_1+a_0)/3$ by (\ref{dPcondition}).
\end{itemize}

If $a_2=(4a_1+a_0)/3$ or $a_2=5a_1/4+a_0/2$, then $a_4<420$ by Theorem~$QS$ and Lemma~$2$. If $a_2=3a_1/4+a_0/2>a_1$ or $a_2=2(a_1+a_0)/3>a_1$, then $a_1<2a_0$. Then $a_4<144$ by Theorem~$QS$ and Lemma~$2$.\\

If $a_2=a_1+a_0/3$, then $d_1=3a_1$, $d_2=4a_1-2a_0/3$. By Theorem~$QS$~$(1)$, $i=0$, either $12a_1\divby a_0$ or $9a_1\divby a_0$. Since $\gcd(a_0/3,a_1)=1$ by Theorem~$WF$, $a_0/3\in \{ 1,2,3,4 \}$. Note that if $6\mid a_0$, then the condition of Theorem~$QS$~$(1)$, $i=0$, is violated. Hence \begin{itemize}
\item either $a_0=3$, $a_2=a_1+1$, $a_3=2a_1-1$, $a_4=3a_1-3$, $d_1=3a_1$, $d_2=4a_1-2$, where $a_1\equiv 1\pmod 3$,
\item or $a_0=9$, $a_2=a_1+3$, $a_3=2a_1-3$, $a_4=3a_1-9$, $d_1=3a_1$, $d_2=4a_1-6$, where $a_1\equiv -1\pmod 3$.
\end{itemize}

These solutions appear in Table~\ref{table:1} (No.~$29$ and No.~$44$).\\

If $a_2=a_1+2a_0/3$, then $a_4=3a_1$, $d_1=3a_1+a_0$, $d_2=4a_1+2a_0/3$. By Theorem~$WF$, $2a_0\divby a_1$, i.e. $a_1=2a_0$, $a_2=8a_0/3$, and so $a_4=6a_0\leq 18$ by Lemma~$2$.\\

If $a_2=a_1+a_0/2$, then $d_1=3a_1+a_0/2$, $d_2=4a_1$. By Theorem~$QS$~$(1)$, $i=0$, either $6a_1\divby a_0$ or $4a_1\divby a_0$. Since $\gcd(a_0/2,a_1)=1$ by Theorem~$WF$, $a_0/2\in \{ 1,2,3 \}$. Note that if $a_0/2$ is odd, then a condition of Theorem~$WF$ is violated. Hence 
$$
a_0=4,\; a_1=2t+1,\; a_2=2t+3,\; a_3=4t+2,\; a_4=6t+1,\; d_1=6t+5,\; d_2=8t+4.
$$

This solution appears in Table~\ref{table:1} (No.~$33$).\\

If $a_3=2a_2-a_0/2$, then $d_1=3a_2$, $d_2=4a_2-a_0$, $a_2<a_1+a_0/2$. By Theorem~$QS$~$(1)$, $i=1$, 
\begin{itemize}
\item either $3a_2\divby a_1$, i.e. $a_2=4a_1/3$ by (\ref{dPcondition}),
\item or $4a_2-a_0\divby a_1$, i.e. $a_2=a_1+a_0/4$ by (\ref{dPcondition}),
\item or $3a_2-a_0\divby a_1$, i.e. $a_2=a_1+a_0/3$ by (\ref{dPcondition}).
\end{itemize}

If $a_2=4a_1/3$, then $d_1=4a_1$, $d_2=16a_1/3-a_0$, $a_1<3a_0/2$. Then $a_4<240$ by Theorem~$QS$ and Lemma~$2$.\\

If $a_2=a_1+a_0/3$, then $a_4=3a_1$, $d_1=3a_1+a_0$, $d_2=4a_1+a_0/3$. By Theorem~$WF$, $a_0\divby a_1$, i.e. $a_0=a_1$, which violates (\ref{AllNonequal}).\\

If $a_2=a_1+a_0/4$, then $d_1=3a_1+3a_0/4$, $d_2=4a_1$. By Theorem~$QS$~$(1)$, $i=0$, $12a_1\divby a_0$, i.e. $a_0=4$ or $a_0=12$. Hence \begin{itemize}
\item either $a_0=4$, $a_2=a_1+1$, $a_3=2a_1$, $a_4=3a_1-1$, $d_1=3a_1+3$, $d_2=4a_1$,
\item or $a_0=12$, $a_2=a_1+3$, $a_3=2a_1$, $a_4=3a_1-3$, $d_1=3a_1+9$, $d_2=4a_1$.
\end{itemize}

Both these items violate conditions of Theorem~$WF$.\\

\subsubsection{Subcase $a_3=a_2-a_0$, $a_4=a_2+a_1-2a_0$}

In this subcase, $d_1=a_2+a_1-a_0$, $d_2=2a_2+a_1-2a_0$. By Theorem~$QS$~$(1)$, $i=2$, $a_1=2a_0$, and hence $a_4=a_2$, which contradicts to (\ref{AllNonequal}).\\

\subsubsection{Subcase $a_3=a_2+a_1-2a_0$, $a_4=a_2+2a_1-3a_0$}

In this subcase, $d_1=a_2+2a_1-2a_0$, $d_2=2a_2+2a_1-3a_0$. By Theorem~$QS$~$(1)$, $i=2$, 
\begin{itemize}
\item either $2(a_1-a_0)\divby a_2$, i.e. $a_2=2(a_1-a_0)$,
\item or $2a_1-3a_0\divby a_2$, i.e. $a_2=2a_1-3a_0$.
\end{itemize}

In both cases, $a_4<40$ by Theorem~$QS$~$(1)$, $i=1$, and Lemma~$2$.\\

\subsubsection{Subcase $a_3=2a_2-a_1-a_0$, $a_4=3a_2-a_1-2a_0$}

In this subcase, $d_1=3a_2-a_1-a_0$, $d_2=4a_2-a_1-2a_0$, $a_2<a_1+a_0$. By Theorem~$QS$~$(1)$, $i=1$, 
\begin{itemize}
\item either $3a_2-a_0\divby a_1$, i.e. $a_2=(4a_1+a_0)/3$ by (\ref{dPcondition}),
\item or $4a_2-2a_0\divby a_1$, i.e. $a_2=5a_1/4+a_0/2$ by (\ref{dPcondition}),
\item or $3a_2-2a_0\divby a_1$, i.e. $a_2=a_1+2a_0/3$ by (\ref{dPcondition}).
\end{itemize}

If $a_2=(4a_1+a_0)/3$ or $a_2=5a_1/4+a_0/2$, then $a_4<220$ by Theorem~$QS$ and Lemma~$2$.\\

If $a_2=a_1+2a_0/3$, then $a_4=2a_1$, $d_1=2a_1+a_0$, $d_2=3a_1+2a_0/3$. By Theorem~$WF$, $2a_0\divby a_1$, i.e. $a_1=2a_0$, $a_2=8a_0/3$, $a_4=4a_0\leq 12$ by Lemma~$2$.\\

\subsection{Case $(12)$}

\subsubsection{Subcase $a_4=2a_3-a_1$}

In this subcase, $d_1=2a_3$, $d_2=2a_3+a_2-a_1$, $a_3<a_0+a_1$. By Theorem~$QS$~$(1)$, $i=2$, 
\begin{itemize}
\item either $2a_3\divby a_2$, i.e. $a_3=3a_2/2$ by (\ref{dPcondition}),
\item or $2a_3-a_1\divby a_2$, i.e. $a_3=a_2+a_1/2$ by (\ref{dPcondition}),
\item or $2a_3-a_0\divby a_2$, i.e. $a_3=a_2+a_0/2$ by (\ref{dPcondition}).
\end{itemize}

In each case, $a_4<135$ by Theorem~$QS$, Theorem~$WF$ and Lemma~$2$.\\

\subsubsection{Subcase $a_4=2a_3-a_2$}

In this subcase, $d_1=2a_3+a_1-a_2$, $d_2=2a_3$, $a_3<a_0+a_2$. By Theorem~$QS$~$(1)$, $i=2$, 
\begin{itemize}
\item either $2a_3+a_1\divby a_2$, i.e. $a_3=(3a_2-a_1)/2$ or $a_3=2a_2-a_1/2$ by (\ref{dPcondition}),
\item or $2a_3\divby a_2$, i.e. $a_3=3a_2/2$ by (\ref{dPcondition}),
\item or $a_3+a_1\divby a_2$, i.e. $a_3=2a_2-a_1$ by (\ref{dPcondition}),
\item or $2a_3+a_1-a_0\divby a_2$, i.e. $a_3=(3a_2+a_0-a_1)/2$ by (\ref{dPcondition}).
\end{itemize}

If $a_3=3a_2/2$ or $a_3=2a_2-a_1/2$ or $a_3=(3a_2+a_0-a_1)/2$, then $a_4<168$ by Theorem~$QS$ and Lemma\,$2$.\\

If $a_3=2a_2-a_1$, then $d_1=3a_2-a_1$, $d_2=4a_2-2a_1$, $a_2<a_0+a_1$. By Theorem~$QS$~$(1)$, $i=1$, 
\begin{itemize}
\item either $3a_2\divby a_1$, i.e. $a_2=\lambda a_1/3$, $\lambda \in \{ 4,5 \}$ by (\ref{dPcondition}),
\item or $4a_2\divby a_1$, i.e. $a_2=\lambda a_1/4$, $\lambda \in \{ 5,6,7 \}$ by (\ref{dPcondition}),
\item or $4a_2-a_0\divby a_1$, i.e. $a_2=(\lambda a_1+a_0)/4$, $\lambda \in \{ 4,5,6 \}$ by (\ref{dPcondition}).
\end{itemize}

Theorem~$QS$ and Lemma~$2$ give bound $a_4<264$ in all cases except for
$$
a_2=a_1+a_0/4.
$$

In this case, $d_1=2a_1+3a_0/4$, $d_2=2a_1+a_0$. By Theorem~$QS$~$(1)$, $i=0$, $8a_1\divby a_0$. Since $\gcd(a_0/4,a_1)=1$ by Theorem~$WF$, $a_0=4$ or $a_0=8$. Hence \begin{itemize}
\item either $a_0=4$, $a_2=a_1+1$, $a_3=a_1+2$, $a_4=a_1+3$, $d_1=2a_1+3$, $d_2=2a_1+4$,
\item or $a_0=8$, $a_2=a_1+2$, $a_3=a_1+4$, $a_4=a_1+6$, $d_1=2a_1+6$, $d_2=2a_1+8$, where $a_1$ is odd.
\end{itemize}

If $a_0=4$, then a condition of Theorem~$WF$ is violated. Hence
$$
a_0=8,\; a_1=2t+1,\; a_2=2t+3,\; a_3=2t+5,\; a_4=2t+7,\; d_1=4t+8,\; d_2=4t+10,\; \mbox{ where $t$ is even}.
$$

This solution appears in Table~\ref{table:1} (No.~$42$).\\

If $a_3=(3a_2-a_1)/2$, then $d_1=2a_2$, $d_2=3a_2-a_1$, $a_2<a_1+2a_0$. By Theorem~$QS$~$(1)$, $i=1$, 
\begin{itemize}
\item either $2a_2\divby a_1$, i.e. $a_2=\lambda a_1/2$, $\lambda \in \{ 3,4,5 \}$ by (\ref{dPcondition}),
\item or $3a_2\divby a_1$, i.e. $a_2=\lambda a_1/3$, $4\leq \lambda \leq 8$ by (\ref{dPcondition}),
\item or $3a_2-a_0\divby a_1$, i.e. $a_2=(\lambda a_1+a_0)/3$, $3\leq \lambda \leq 7$ by (\ref{dPcondition}).
\end{itemize}

Theorem~$QS$ and Lemma~$2$ give bound $a_4<221$ in all cases except for
$$
a_2=a_1+a_0/3.
$$

In this case, $d_1=2a_1+2a_0/3$, $d_2=2a_1+a_0$. By Theorem~$QS$~$(1)$, $i=0$, $6a_1\divby a_0$, i.e. $a_0=6$, because $6\mid a_0$ and $\gcd(a_0/3,a_1)=1$ by Theorem~$WF$. Hence
$$
a_0=6,\; a_2=a_1+2,\; a_3=a_1+3,\; a_4=a_1+4,\; d_1=2a_1+4,\; d_2=2a_1+6,\; \mbox{ where }\; a_1\equiv 1 \pmod{3} \; \mbox{ is odd}.
$$ 

This solution appears in Table~\ref{table:1} (No.~$38$).\\

\subsubsection{Subcase $a_2=2a_1-a_0$, $a_4=2a_3+a_0-a_1$}

In this subcase, $a_4<a_0+a_3$ by (\ref{dPcondition}), which implies that $a_3<a_1$. This contradicts to (\ref{AllNonequal}).\\  

\subsubsection{Subcase $a_3=2a_2-a_1-a_0$, $a_4=3a_2-2a_1-a_0$}

In this subcase, $a_4<a_0+a_3$ by (\ref{dPcondition}), which implies that $a_2<a_0+a_1$. Hence $a_3<a_2$, which contradicts to (\ref{AllNonequal}).\\ 

This finishes the proof of the Main Theorem.\\

\section*{Acknowledgement}

We are grateful to Beijing International Center for Mathematical Research, Simons Foundation and Peking University for support, excellent working conditions and encouraging atmosphere.\\

\newpage
\section*{Appendix}

\begin{center}
\begin{longtable}{@{\extracolsep{\fill}} c c c l @{}}
\caption{Infinite series of well-formed quasi-smooth codimension $2$ weighted complete intersection del Pezzo surfaces, which are not intersections with linear cones.}\label{table:1}\\
\toprule[.1em]
No. & $(a_0,\; a_1,\; a_2,\; a_3,\; a_4; \; d_1,\; d_2)$ & $I$ & conditions on the parameters\\ [1ex]
\midrule[.1em]
$1$ & $\begin{array}{c} (a_0,\; a_1,\;(\nu a_0 -1)a_1,\\ \nu a_0a_1-\frac{a_0+a_1}{2},\; (\nu a_1-1)a_0; \\ \nu a_0a_1,\; 2\nu a_0a_1 - a_0-a_1)\end{array}$  & $\dfrac{a_0+a_1}{2}$ & $\begin{array}{c}\gcd\left( \nu a_0-1,\frac{a_1-a_0}{2}\right)=1, \\ \gcd(a_0,a_1)=1,\; a_0, a_1 \; \mbox{odd} \\ \nu\geq \max(1,3-a_0) \end{array}$\\\hline
$2$ & $\begin{array}{c} (a_0,\; a_1,\; \left( \frac{\nu a_0 -1}{2}\right) a_1,\\ \left( \frac{\nu a_1-1}{2}\right) a_0, \; \left(\frac{\nu a_0+1}{2}\right) a_1-a_0; \\ \left( \frac{\nu a_0+1}{2}\right) a_1,\; (\nu a_1-1)a_0)\end{array}$ & $\dfrac{a_0+a_1}{2}$ & $\begin{array}{c} \gcd\left(\frac{\nu a_0-1}{2},\frac{a_1-a_0}{2}\right)=1,\\ \gcd(a_0,a_1)=1,\; \nu, a_0, a_1\; \mbox{odd}\\ \nu\geq \max(1,4-a_0) \end{array}$\\\hline
$3$ & $\begin{array}{c} (2b_0,\; 2b_1,\; (\nu b_0 -1)b_1+b_0, \\ \nu b_0b_1, \; (\nu b_0+1)b_1-b_0; \\ (\nu b_0+1)b_1+b_0,\; 2\nu b_0b_1)\end{array}$ & $b_0+b_1$ & $\begin{array}{c} \gcd(\nu ,b_1-b_0)=1,\\ \gcd(b_0,b_1)=1,\; \nu, b_0, b_1\; \mbox{odd} \\ \nu\geq \max(1,4-b_0) \end{array}$\\\hline
$4$ & $\begin{array}{c} (a_0,\; a_1,\; \nu a_0a_1+\frac{a_0-a_1}{2},\\ \nu a_0a_1,\; \nu a_0a_1+\frac{a_1-a_0}{2}; \\ \nu a_0a_1 +\frac{a_0+a_1}{2},\; 2\nu a_0a_1)\end{array}$ & $\dfrac{a_0+a_1}{2}$ & $\begin{array}{c} \gcd\left(\nu,\frac{a_1-a_0}{2}\right)=1,\\ \gcd(a_0,a_1)=1,\; a_0, a_1\; \mbox{odd}  \\ \nu\geq \max(1,3-a_0)\end{array}$\\\hline
$5$ & $\begin{array}{c} (a_0,\; a_1,\; \left( \frac{\nu a_1+1}{2}\right) a_0-a_1,\\ \left( \frac{\nu a_0-1}{2}\right) a_1,\; \left( \frac{\nu a_1-1}{2}\right) a_0; \\ \left( \frac{\nu a_1+1}{2}\right) a_0,\; (\nu a_0-1)a_1)\end{array}$ & $\dfrac{a_0+a_1}{2}$ & $\begin{array}{c} \gcd\left(\frac{\nu a_0-1}{2},\frac{a_1-a_0}{2}\right)=1,\\ \gcd(a_0,a_1)=1,\; \nu, a_0, a_1\; \mbox{odd}\\ \nu\geq \max(1,6-a_0) \end{array}$\\\hline
$6$ & $\begin{array}{c} (a_0,\; a_1,\; 2a_1-a_0,\\ (\nu a_0-1)(2a_1-a_0),\; \nu a_0 (2a_1-a_0)-a_1; \\ (\nu a_0-1)(2a_1-a_0)+a_1,\; \nu a_0 (2a_1-a_0))\end{array}$ & $a_1$ & $\begin{array}{c} \gcd(\nu a_0-1,a_1-a_0)=1,\\ \gcd(a_0,a_1)=1,\; a_0\; \mbox{odd} \\ \nu\geq \max(1,3-a_0) \end{array}$\\\hline
$7$ & $\begin{array}{c} (a_0,\; a_1,\; 2a_1-a_0,\; \left(\frac{\nu a_0-1}{2}\right)(2a_1-a_0)-a_1,\\ \left(\frac{\nu a_0-1}{2}\right)(2a_1-a_0)-a_0; \\ \left(\frac{\nu a_0-1}{2}\right)(2a_1-a_0),\; \left(\frac{\nu a_0+1}{2}\right)(2a_1-a_0)-a_1)\end{array}$ & $a_1$ & $\begin{array}{c} \gcd\left(\frac{\nu a_0-3}{2},a_1-a_0\right)=1,\\ \gcd(a_0,a_1)=1,\; \nu, a_0\; \mbox{odd}\\ \nu\geq \max(1,6-a_0) \end{array}$\\\hline
$8$ & $\begin{array}{c} (a_0,\; a_1,\; 2a_1-a_0,\; \left(\frac{\nu a_0-1}{2}\right)(2a_1-a_0),\\ \left(\frac{\nu a_0-1}{2}\right)(2a_1-a_0)+a_1-a_0; \\ \left(\frac{\nu a_0-1}{2}\right)(2a_1-a_0)+a_1,\; \left(\frac{\nu a_0+1}{2}\right)(2a_1-a_0))\end{array}$ & $a_1$ & $\begin{array}{c} \gcd\left(\frac{\nu a_0-1}{2},a_1-a_0\right)=1,\\ \gcd(a_0,a_1)=1,\; \nu, a_0\; \mbox{odd}\\ \nu\geq \max(1,4-a_0) \end{array}$\\\hline
$9$ & $\begin{array}{c} (a_0,\; a_1,\; 2a_1-a_0,\\ \nu a_0(2a_1-a_0)-a_1,\; \nu a_0 (2a_1-a_0)-a_0; \\ \nu a_0(2a_1-a_0),\; (\nu a_0+1) (2a_1-a_0)-a_1)\end{array}$ & $a_1$ & $\begin{array}{c} \gcd(\nu a_0-1,a_1-a_0)=1,\\ \gcd(a_0,a_1)=1,\; a_0\; \mbox{odd}\\ \nu\geq \max(1,3-a_0) \end{array}$\\\hline
$10$ & $\begin{array}{c} (2b_0,\; 2b_1,\; (2\nu b_0 -1)b_1,\\ (2\nu b_1-1)b_0,\; 2(2\nu b_0b_1- b_0-b_1); \\ 2(2\nu b_0-1)b_1,\; 2(2\nu b_1-1)b_0)\end{array}$ & $b_0+b_1$ & $\begin{array}{c} \gcd(b_0,b_1)=1,\\ b_0, b_1\; \mbox{odd}\\ \nu\geq \max(1,3-b_0) \end{array}$\\\hline
$11$ & $\begin{array}{c} (a_0,\; a_1,\; \left(\frac{\nu a_0 -1}{2}\right)a_1, \\ \left(\frac{\nu a_1 -1}{2}\right)a_0,\; \nu a_0a_1- a_0-a_1; \\ (\nu a_0-1)a_1,\; (\nu a_1-1)a_0)\end{array}$ & $\dfrac{a_0+a_1}{2}$ & $\begin{array}{c} \gcd(a_0,a_1)=1,\\ \nu, a_0, a_1\; \mbox{odd}\\ \nu\geq \max(1,4-a_0)\end{array}$\\\hline
$12$ & $\begin{array}{c} (a_0,\; a_1,\; (\nu a_0-1)a_1, \\ \nu a_0a_1-\frac{a_0+a_1}{2},\; 2\nu a_0a_1- a_0-2a_1; \\ 2(\nu a_0-1)a_1,\; 2\nu a_0a_1-a_0-a_1)\end{array}$ & $\dfrac{a_0+a_1}{2}$ & $\begin{array}{c} \gcd(a_0,a_1)=1,\\ a_0, a_1\; \mbox{odd}\\ \nu\geq \max(1,3-a_0) \end{array}$\\\hline
$13$ & $\begin{array}{c} (a_0,\; a_1,\; (\nu a_0-1)a_1+\frac{a_1-a_0}{2}, \\ (\nu a_1-1)a_0,\; 2\nu a_0a_1- 2a_0-a_1; \\ 2\nu a_0a_1-a_0-a_1,\; 2(\nu a_1-1)a_0)\end{array}$ & $\dfrac{a_0+a_1}{2}$ & $\begin{array}{c} \gcd(a_0,a_1)=1,\\ a_0, a_1\; \mbox{odd}\\ \nu\geq \max(1,3-a_0)\end{array}$\\\hline
$14$ & $\begin{array}{c} (a_0,\; a_1,\; \left(\frac{\nu a_1 -1}{2}\right)a_0-a_1, \\ \left(\frac{\nu a_0 -1}{2}\right)a_1-a_0,\; \nu a_0a_1- 2(a_0+a_1); \\ \nu a_0a_1-a_0-2a_1,\; \nu a_0a_1-2a_0-a_1)\end{array}$ & $\dfrac{a_0+a_1}{2}$ & $\begin{array}{c} \gcd(a_0,a_1)=1,\\ \nu, a_0, a_1\; \mbox{odd}\\ \nu\geq \max(1,6-a_0)\end{array}$\\\hline
$15$ & $(1,\; 1,\; t,\; t,\; 2t-1;\; 2t,\; 2t)$ & $1$ & \\\hline
$16$ & $(1,\; 2,\; 2t+1,\; 2t+1,\; 4t+1;\; 4t+2,\; 4t+3)$ & $1$ & \\\hline
$17$ & $(1,\; t,\; 2t-1,\; 2t-1,\; 3t-2;\; 3t-1,\; 4t-2)$ & $t$ & \\\hline
$18$ & $(1,\; 4t-2,\; 6t-3,\; 9t-5,\; 12t-7;\; 12t-6,\; 18t-10)$ & $t$ & \\\hline
$19$ & $(1,\; 2t-1,\; 2t-1,\; 3t-2,\; 4t-3;\; 4t-2,\; 6t-4)$ & $t$ & \\\hline
$20$ & $(1,\; 6t-1,\; 8t-2,\; 12t-3,\; 18t-5;\; 18t-4,\; 24t-6)$ & $2t$ & \\\hline
$21$ & $(2,\; 2,\; 2t+1,\; 2t+1,\; 4t;\; 4t+2,\; 4t+2)$ & $2$ & \\\hline
$22$ & $(2,\; 3,\; 3t,\; 3t+1,\; 3t+1;\; 3t+3,\; 6t+2)$ & $2$ & \\\hline
$23$ & $(2,\; 3,\; 3t+1,\; 3t+2,\; 6t+1;\; 6t+3,\; 6t+4)$ & $2$ & \\\hline
$24$ & $(2,\; 4,\; 2t+3,\; 2t+3,\; 4t+4;\; 4t+6,\; 4t+8)$ & $2$ & \\\hline
$25$ & $(2,\; 2t+1,\; 2t+1,\; 4t+1,\; 6t+1;\; 6t+3,\; 8t+2)$ & $1$ & \\\hline
$26$ & $(3,\; t+2,\; 2t+1,\; 2t+1,\; 3t;\; 3t+3,\; 4t+2)$ & $t+2$ & $t\not\equiv 1\pmod 3$ \\\hline
$27$ & $(3,\; 3t,\; 3t+1,\; 3t+1,\; 3t+2;\; 6t+2,\; 6t+3)$ & $2$ & \\\hline
$28$ & $(3,\; t+2,\; t+3,\; t+3,\; 2t+3;\; 2t+6,\; 3t+6)$ & $2$ & $3\nmid t$ \\\hline
$29$ & $(3,\; 3t+1,\; 3t+2,\; 6t+1,\; 9t;\; 9t+3,\; 12t+2)$ & $2$ & \\\hline
$30$ & $(3,\; 2t+1,\; 2t+1,\; 3t,\; 4t-1;\; 4t+2,\; 6t)$ & $t+2$ & $t\not\equiv 1\pmod 3$ \\\hline
$31$ & $(4,\; 6,\; 6t+5,\; 6t+7,\; 12t+8;\; 12t+12,\; 12t+14)$ & $4$ & \\\hline
$32$ & $(4,\; 6,\; 6t+3,\; 6t+5,\; 6t+5;\; 6t+9,\; 12t+10)$ & $4$ & \\\hline
$33$ & $(4,\; 2t+3,\; 2t+5,\; 4t+6,\; 6t+7;\; 6t+11,\; 8t+12)$ & $2$ & \\\hline
$34$ & $(4,\; 2t+3,\; 2t+3,\; 2t+5,\; 4t+4;\; 4t+8,\; 6t+9)$ & $2$ & \\\hline
$35$ & $(4,\; 2t+3,\; 2t+3,\; 4t+4,\; 6t+5;\; 6t+9,\; 8t+8)$ & $2$ & \\\hline
$36$ & $(4,\; 2t+3,\; 4t+6,\; 6t+7,\; 8t+8;\; 8t+12,\; 12t+14)$ & $2$ & \\\hline
$37$ & $(4,\; 4t+1,\; 4t+2,\; 4t+3,\; 4t+3;\; 8t+4,\; 8t+6)$ & $3$ & \\\hline
$38$ & $(6,\; 6t+1,\; 6t+3,\; 6t+4,\; 6t+5;\; 12t+6,\; 12t+8)$ & $5$ & \\\hline
$39$ & $(6,\; 6t+3,\; 6t+5,\; 6t+5,\; 6t+7;\; 12t+10,\; 12t+12)$ & $4$ & \\\hline
$40$ & $(6,\; 6t+3,\; 6t+5,\; 6t+5,\; 12t+4;\; 12t+10,\; 18t+9)$ & $4$ & \\\hline
$41$ & $(7,\; 4t+6,\; 6t+9,\; 9t+10,\; 12t+11;\; 12t+18,\; 18t+20)$ & $t+5$ & $t\not\equiv 2 \pmod 7$\\\hline
$42$ & $(8,\; 4t+5,\; 4t+7,\; 4t+9,\; 4t+11;\; 8t+16,\; 8t+18)$ & $6$ & \\\hline
$43$ & $(8,\; 4t+5,\; 4t+7,\; 4t+9,\; 8t+6;\; 8t+14,\; 12t+15)$ & $6$ & \\\hline
$44$ & $(9,\; 3t+8,\; 3t+11,\; 6t+13,\; 9t+15;\; 9t+24,\; 12t+26)$ & $6$ & \\\hline
$45$ & $(9,\; 3t+8,\; 3t+11,\; 3t+14,\; 6t+13;\; 6t+22,\; 9t+27)$ & $6$ & \\
\bottomrule[.1em]
\end{longtable}
\end{center}

In Table~\ref{table:1}, $a_0<a_1$, $b_0<b_1$, $\nu$, $t$ denote positive integers. Note that, when $t\geq 2$, No.~$19$ is a special case of No.~$1$ (when $a_0=1$, $\nu=2$), of No.~$2$ (when $a_0=1$, $\nu=3$), of No.~$11$ (when $a_0=1$, $\nu=3$) and of No.~$12$ (when $a_0=1$, $\nu=2$). No.~$30$ is a special case of No.~$2$ (when $a_0=3$, $\nu=1$) and of No.~$11$ (when $a_0=3$, $\nu=1$). No.~$17$ is a special case of No.~$6$ (when $a_0=1$, $\nu=2$) and of No.~$8$ (when $a_0=1$, $\nu=3$). No.~$26$ is a special case of No.~$8$ (when $a_0=3$, $\nu=1$).\\

Table~\ref{table:2} below lists well-formed quasi-smooth codimension $2$ weighted complete intersection del Pezzo surfaces, which are not intersections with linear cones, with $a_4\leq 500$, $d_2\leq 1000$, which do not appear in Table~\ref{table:1}. It was obtained by running a computer program, whose code is available from the author upon request. According to the Main Theorem, Table~\ref{table:1} and Table~\ref{table:2} together list all well-formed quasi-smooth codimension $2$ weighted complete intersection del Pezzo surfaces, which are not intersections with linear cones.\\

\begin{center}
\begin{longtable}{ c || c || c }
\caption{Sporadic well-formed quasi-smooth codimension $2$ weighted complete intersection del Pezzo surfaces, which are not intersections with linear cones.}\label{table:2}\\
\toprule[.1em]
$(a_0,\; a_1,\; a_2,\; a_3,\; a_4; \; d_1,\; d_2)$ & $(a_0,\; a_1,\; a_2,\; a_3,\; a_4; \; d_1,\; d_2)$ & $(a_0,\; a_1,\; a_2,\; a_3,\; a_4; \; d_1,\; d_2)$\\ [1ex]
\midrule[.1em]
$(13,\; 22,\; 55,\; 76,\; 97;\; 110,\; 152)$ & $(9,\; 19,\; 24,\; 31,\; 53;\; 62,\; 72)$ & $(1,\; 7,\; 11,\; 17,\; 27;\; 28,\; 34)$ \\\hline
$(11,\; 27,\; 36,\; 62,\; 97;\; 108,\; 124)$ & $(14,\; 19,\; 25,\; 32,\; 45;\; 64,\; 70)$ & $(5,\; 7,\; 10,\; 14,\; 23;\; 28,\; 30)$ \\\hline
$(13,\; 18,\; 45,\; 61,\; 77;\; 90,\; 122)$ & $(10,\; 17,\; 25,\; 34,\; 43;\; 60,\; 68)$ & $(2,\; 7,\; 10,\; 13,\; 18;\; 20,\; 28)$ \\\hline
$(11,\; 29,\; 39,\; 49,\; 67;\; 78,\; 116)$ & $(11,\; 17,\; 24,\; 31,\; 37;\; 48,\; 68)$ & $(6,\; 7,\; 9,\; 11,\; 14;\; 18,\; 28)$ \\\hline
$(11,\; 29,\; 38,\; 48,\; 85;\; 96,\; 114)$ & $(11,\; 14,\; 21,\; 33,\; 52;\; 63,\; 66)$ & $(5,\; 8,\; 9,\; 12,\; 19;\; 24,\; 27)$ \\\hline
$(13,\; 23,\; 35,\; 57,\; 79;\; 92,\; 114)$ & $(13,\; 14,\; 23,\; 33,\; 43;\; 56,\; 66)$ & $(2,\; 7,\; 8,\; 13,\; 19;\; 21,\; 26)$ \\\hline
$(13,\; 23,\; 34,\; 56,\; 89;\; 102,\; 112)$ & $(13,\; 14,\; 23,\; 32,\; 33;\; 46,\; 65)$ & $(1,\; 5,\; 9,\; 13,\; 17;\; 18,\; 26)$ \\\hline
$(13,\; 23,\; 35,\; 47,\; 57;\; 70,\; 104)$ & $(11,\; 15,\; 20,\; 32,\; 49;\; 60,\; 64)$ & $(5,\; 6,\; 9,\; 13,\; 13;\; 18,\; 26)$ \\\hline
$(11,\; 25,\; 34,\; 43,\; 57;\; 68,\; 100)$ & $(11,\; 17,\; 24,\; 31,\; 38;\; 55,\; 62)$ & $(6,\; 8,\; 9,\; 11,\; 13;\; 22,\; 24)$ \\\hline
$(11,\; 29,\; 39,\; 49,\; 59;\; 88,\; 98)$ & $(10,\; 13,\; 25,\; 31,\; 37;\; 50,\; 62)$ & $(1,\; 5,\; 8,\; 12,\; 19;\; 20,\; 24)$ \\\hline
$(13,\; 20,\; 31,\; 49,\; 67;\; 80,\; 98)$ & $(11,\; 17,\; 20,\; 27,\; 43;\; 54,\; 60)$ & $(5,\; 7,\; 8,\; 11,\; 14;\; 21,\; 22)$ \\\hline
$(11,\; 25,\; 32,\; 41,\; 71;\; 82,\; 96)$ & $(9,\; 15,\; 23,\; 23,\; 37;\; 46,\; 60)$ & $(2,\; 5,\; 8,\; 11,\; 14;\; 16,\; 22)$ \\\hline
$(13,\; 20,\; 29,\; 47,\; 74;\; 87,\; 94)$ & $(13,\; 14,\; 19,\; 29,\; 44;\; 57,\; 58)$ & $(3,\; 7,\; 8,\; 9,\; 13;\; 16,\; 21)$ \\\hline
$(11,\; 21,\; 28,\; 47,\; 73;\; 84,\; 94)$ & $(14,\; 15,\; 19,\; 26,\; 37;\; 52,\; 56)$ & $(4,\; 5,\; 7,\; 10,\; 13;\; 18,\; 20)$ \\\hline
$(13,\; 14,\; 35,\; 46,\; 57;\; 70,\; 92)$ & $(9,\; 15,\; 23,\; 23,\; 31;\; 46,\; 54)$ & $(1,\; 4,\; 7,\; 10,\; 13;\; 14,\; 20)$ \\\hline
$(13,\; 20,\; 31,\; 42,\; 49;\; 62,\; 91)$ & $(13,\; 14,\; 19,\; 23,\; 29;\; 42,\; 52)$ & $(3,\; 5,\; 6,\; 8,\; 13;\; 16,\; 18)$ \\\hline
$(9,\; 23,\; 30,\; 38,\; 67;\; 76,\; 90)$ & $(11,\; 13,\; 19,\; 25,\; 27;\; 38,\; 52)$ & $(3,\; 5,\; 7,\; 9,\; 11;\; 16,\; 18)$ \\\hline
$(11,\; 18,\; 27,\; 44,\; 61;\; 72,\; 88)$ & $(11,\; 17,\; 20,\; 24,\; 27;\; 44,\; 51)$ & $(2,\; 5,\; 6,\; 9,\; 13;\; 15,\; 18)$ \\\hline
$(11,\; 25,\; 34,\; 43,\; 52;\; 77,\; 86)$ & $(5,\; 14,\; 17,\; 21,\; 37;\; 42,\; 51)$ & $(3,\; 5,\; 7,\; 9,\; 11;\; 14,\; 18)$ \\\hline
$(10,\; 19,\; 35,\; 43,\; 51;\; 70,\; 86)$ & $(11,\; 13,\; 19,\; 25,\; 31;\; 44,\; 50)$ & $(1,\; 4,\; 5,\; 7,\; 11;\; 12,\; 15)$ \\\hline
$(14,\; 17,\; 29,\; 41,\; 44;\; 58,\; 85)$ & $(9,\; 12,\; 19,\; 19,\; 29;\; 38,\; 48)$ & $(4,\; 5,\; 6,\; 7,\; 8;\; 12,\; 14)$ \\\hline
$(11,\; 21,\; 29,\; 37,\; 47;\; 58,\; 84)$ & $(1,\; 9,\; 15,\; 23,\; 23;\; 24,\; 46)$ & $(3,\; 4,\; 6,\; 7,\; 8;\; 12,\; 14)$ \\\hline
$(14,\; 17,\; 29,\; 41,\; 53;\; 70,\; 82)$ & $(9,\; 12,\; 19,\; 19,\; 26;\; 38,\; 45)$ & $(3,\; 4,\; 5,\; 6,\; 7;\; 11,\; 12)$ \\\hline
$(13,\; 17,\; 27,\; 41,\; 55;\; 68,\; 82)$ & $(9,\; 10,\; 15,\; 22,\; 23;\; 32,\; 45)$ & $(3,\; 4,\; 5,\; 6,\; 7;\; 10,\; 12)$ \\\hline
$(14,\; 17,\; 27,\; 39,\; 64;\; 78,\; 81)$ & $(10,\; 11,\; 15,\; 22,\; 29;\; 40,\; 44)$ & $(2,\; 3,\; 5,\; 6,\; 7;\; 10,\; 12)$ \\\hline
$(11,\; 21,\; 26,\; 34,\; 57;\; 68,\; 78)$ & $(11,\; 13,\; 14,\; 20,\; 29;\; 40,\; 42)$ & $(3,\; 3,\; 5,\; 5,\; 7;\; 10,\; 12)$ \\\hline
$(13,\; 20,\; 29,\; 31,\; 47;\; 60,\; 78)$ & $(5,\; 9,\; 12,\; 20,\; 31;\; 36,\; 40)$ & $(2,\; 3,\; 4,\; 5,\; 5;\; 8,\; 10)$ \\\hline
$(13,\; 17,\; 27,\; 37,\; 41;\; 54,\; 78)$ & $(1,\; 8,\; 13,\; 19,\; 31;\; 32,\; 39)$ & $(1,\; 3,\; 3,\; 5,\; 5;\; 6,\; 10)$ \\\hline
$(13,\; 17,\; 24,\; 38,\; 59;\; 72,\; 76)$ & $(2,\; 9,\; 12,\; 19,\; 19;\; 21,\; 38)$ & $(2,\; 2,\; 3,\; 3,\; 3;\; 6,\; 6)$ \\\hline
$(11,\; 25,\; 32,\; 34,\; 41;\; 66,\; 75)$ & $(9,\; 11,\; 12,\; 17,\; 25;\; 34,\; 36)$ & $(1,\; 2,\; 2,\; 3,\; 3;\; 4,\; 6)$ \\\hline
$(11,\; 21,\; 29,\; 37,\; 45;\; 66,\; 74)$ & $(1,\; 7,\; 12,\; 17,\; 23;\; 24,\; 35)$ & \\
\bottomrule[.1em]
\end{longtable}
\end{center}

\newpage
\bibliographystyle{ams-plain}

\bibliography{delPezzoWCI}

\providecommand{\bysame}{\leavevmode\hbox to3em{\hrulefill}\thinspace}
\begin{thebibliography}{1}

\bibitem{Boyer}
C.~Boyer, K.~Galicki, and M.~Nakamaye, \emph{On the geometry of
  {S}asakian-{E}instein $5$-manifolds}, Mathematische {A}nnalen \textbf{325}
  (2003), no.~3, 485--524.

\bibitem{CheltsovShramov}
I.~Cheltsov and C.~Shramov, \emph{Del {P}ezzo zoo}, Experimental {M}athematics
  \textbf{22} (2013), no.~3, 313--326.

\bibitem{Chen3}
J.-J. Chen, J.~A. Chen, and M.~Chen, \emph{On quasismooth weighted complete
  intersections}, Journal of {A}lgebraic {G}eometry \textbf{20} (2011), no.~2,
  239--262.

\bibitem{Dolgachev}
I.~Dolgachev, \emph{Weighted projective varieties}, Group actions and vector
  fields ({V}ancouver, {B}.{C}., $1981$), Lecture {N}otes in {M}athematics,
  vol. 956, Springer, Berlin, 1982, pp.~34--71.

\bibitem{Fletcher}
A.~R. Iano-Fletcher, \emph{Working with weighted complete intersections},
  {E}xplicit birational geometry of $3$-folds, London {M}athematical {S}ociety
  {L}ecture {N}ote {S}eries, vol. 281, Cambridge {U}niversity {P}ress,
  Cambridge, 2000, pp.~101--173.

\bibitem{JohnsonKollar}
J.~M. Johnson and J.~Koll{\'a}r, \emph{K{\"a}hler-{E}instein metrics on log del
  {P}ezzo surfaces in weighted projective $3$-spaces}, Universit{\'e} de
  {G}renoble. {A}nnales de l'{I}nstitut {F}ourier \textbf{51} (2001), no.~1,
  69--79.

\bibitem{KimPark}
I.-K. Kim and J.~Park, \emph{Log canonical thresholds of complete intersection
  log del {P}ezzo surfaces}, Proceedings of the {E}dinburgh {M}athematical
  {S}ociety. {S}eries {II} \textbf{58} (2015), no.~2, 445--483.

\end{thebibliography}

\end{document}